\documentclass[twoside]{article}

\usepackage[accepted]{aistats2022}

\usepackage[utf8]{inputenc}

\usepackage[ruled,vlined]{algorithm2e}
\usepackage{amsfonts}
\usepackage{amsmath}
\usepackage{amssymb}
\usepackage{amsthm}
\usepackage{caption}
\usepackage{subcaption}
\usepackage{dsfont}
\usepackage{hyperref}
\hypersetup{
    colorlinks=true,
    linkcolor=blue,
    filecolor=magenta,      
    urlcolor=blue,
}
\usepackage[table]{xcolor}
\usepackage{optidef}
\usepackage{tikz}
\usetikzlibrary{shapes, arrows.meta, positioning, matrix, fit, tikzmark}
\usepackage{verbatim}
\usepackage{float}
\usepackage{rotating}
\usepackage{multicol}
\usepackage{tabularx}
\usepackage{nccmath}
\usepackage[round]{natbib}

\colorlet{myRed}{red!20}
 
\newcommand{\quadmat}[3]{\scriptsize \begin{pmatrix} #1 & #2 \\ #2 ^\top & #3 \end{pmatrix}}
\newcommand{\innerMat}[2]{\left \langle #1 , \, #2 \right \rangle}

\usepackage{aligned-overset}

\newif\ifoutlinenotes %
\outlinenotestrue

\ifoutlinenotes

\else

\fi

\newtheorem{theorem}{Theorem}[section]
\newtheorem*{theorem*}{Theorem}
\newtheorem*{corollary*}{Corollary}
\newtheorem{corollary}{Corollary}[theorem]

\newcommand{\getdim}[1]%
{ \path (#1.south west);
  \pgfgetlastxy{\xsw}{\ysw}
  \path (#1.north east);
  \pgfgetlastxy{\xne}{\yne}
  \pgfmathsetlengthmacro{\nodewidth}{\xne-\xsw}
  \pgfmathsetlengthmacro{\nodeheight}{\yne-\ysw}

}

\begin{document}
\twocolumn[
\aistatstitle{A Unified View of SDP-based Neural Network Verification through Completely Positive Programming}
\aistatsauthor{ Robin Brown \And Edward Schmerling \And  Navid Azizan \And Marco Pavone }
\aistatsaddress{Stanford University \\ \texttt{rabrown1@stanford.edu}
\And  Stanford University \\ \texttt{schmrlng@stanford.edu}
\And MIT \\ \texttt{azizan@mit.edu}
\And Stanford University \\ \texttt{pavone@stanford.edu} } ]

\begin{abstract}

Verifying that input-output relationships of a neural network conform to prescribed operational specifications is a key enabler towards deploying these networks in safety-critical applications.
Semidefinite programming (SDP)-based approaches to Rectified Linear Unit (ReLU) network verification transcribe this problem into an optimization problem, where the accuracy of any such formulation reflects the level of fidelity in how the neural network computation is represented, as well as the relaxations of intractable constraints. 
While the literature contains much progress on improving the tightness of SDP formulations while maintaining tractability, comparatively little work has been devoted to the other extreme, i.e., how to most accurately capture the original verification problem before SDP relaxation.
In this work, we develop an exact, convex formulation of verification as a completely positive program (CPP), and provide analysis showing that our formulation is minimal---the removal of any constraint fundamentally misrepresents the neural network computation. 
We leverage our formulation to provide a unifying view of existing approaches, and give insight into the source of large relaxation gaps observed in some cases.

\end{abstract}
\section{INTRODUCTION}

While neural networks today empower many consumer products in image and natural language understanding, they have been shown to fail in surprising and unexpected ways, potentially deterring broad deployment in safety critical settings. 
One avenue for inspiring confidence in neural networks is through the formal verification of safety rules specified as input-output relationships describing limits on the expected behavior of a network.
In this paper, we consider algorithms that pose verification as an optimization problem, where the objective encodes a metric of satisfaction of the safety rule, the constraints represent the neural network computation, and the optimal objective value directly corresponds to confirmation or denial of the safety rule.

Sound and complete, or exact, verifiers must always return the correct decision, which inherently requires exact representation of the neural network computations.
Due to the NP-completeness of verification \cite{KatzBarrettEtAl2017}, exact verifiers face a complexity barrier prohibiting faster than exponential run-time in the worst case.
This suggests that faithfully representing the forward pass of a neural network is at odds with tractable optimization formulations.

Sound but incomplete, or relaxed, verifiers must never return a false assertion of safety, but may conservatively suggest a network is unsafe when it is truly safe.
The conservatism is a result of approximating the neural network computation, and is the trade-off for improved tractability.
In the context of optimization-based verification, this is achieved by loosening exact constraints into approximate constraints; the mismatch between corresponding optimal objective values is termed the relaxation gap.
In this work we aim to develop a deeper understanding of how to tune the balance between tightness and efficiency, motivated by the central challenge of devising a systematic family of relaxations with exact representation of neural network computations as a limiting case.

\paragraph{Contributions.} Our contributions are twofold:
\begin{enumerate}
    \item We develop an \emph{exact, convex formulation} of the verification problem as a completely positive program (CPP); these are linear optimization problems over the cone of completely positive matrices.
    While the complexity of verification is not resolved by the proposed formulation, it is packaged entirely in the complete positivity constraint, with the neural network computation being exactly represented by tractable linear constraints.
    This gives a clean separation between the two competing desiderata of accuracy and tractability in relaxed verification, opening the door for new classes of verifiers that predictably trade-off tightness and efficiency.
    
    \item 
    We also provide analysis explaining how properties of the CPP formulation evolve when the complete positivity constraint is relaxed. 
    We find that many of the favorable properties of the CPP formulation are retained in SDP relaxations, showing that it is a convenient starting point for constructing tight SDP-based verifiers.
    Finally, we contextualize existing work in this shared framework, clearly laying out their similarities and differences with the proposed framework.

\end{enumerate}

\subsection{Related Work}

\paragraph{Sound and Complete Verifiers.}
A number of exact verifiers rely on reformulating verification as a mixed integer linear program (MILP) \citet{ChengNuehrenbergEtAl2017,LomuscioMaganti2017,FischettiJo2017,TjengXiaoEtAl2019} or Satifiability Modulo Theories (SMT) problem \cite{ScheiblerWintererEtAl2015} and calling an off-the-shelf solver. 
Others have developed custom solvers specifically meant to leverage the structure of the verification problem \cite{KatzBarrettEtAl2017,Ehlers2017,DePalmaBunelEtAl2021,JaeckleLuEtAl2021}.
\cite{BunelTurkaslanEtAl2018} provides a unified perspective of exact verifiers through the lens of branch and bound. 
The primary differences across the various methods stem from (1) the way bounds are computed, (2) the type of branching that is considered, and (3) strategies to guide branching. 
Whether deferring to the decisions of an off-the-shelf solver, or employing verification-specific decisions, all exact verifiers ultimately rely on exhaustive search.

\paragraph{Sound but Incomplete Verifiers.} 
Sound but incomplete verifiers bypass exhaustive search by solving convex relaxations of the verification problem, where the neuron values either appear as or can be derived from variables in the optimization problem.
In our discussion, we distinguish between first-order and second-order relaxations. 
Broadly speaking, we categorize a relaxation as first order if the optimization variables represent degree one monomials of the neuron values---in this case, the neuron values can be directly read off the from the optimal solution. 
We consider a relaxation to be second order if the optimization variables represent degree two monomials of the neuron values---in this case, recovering neuron values typically requires factoring the optimal solution. 
Because relaxed verifiers rely on approximating the neural network computations, the optimal neuron values may not necessarily correspond with a valid forward pass of the network.

Under the category of first-order relaxations are methods that are based on linear outer bounds of activations functions \cite{Ehlers2017,SinghGehrEtAl2018}, interval bound propagation \cite{GowalDvijothamEtAl2018,WangPeiEtAl2018,WangPeiEtAl2018b}, and the dual of the relaxed or original non-convex problem \cite{WangPeiEtAl2018b,WongSchmidtEtAl2018,DvijothamStanforthEtAl2018,DvijothamGowalEtAl2018}. 
\cite{SalmanYangEtAl2019} gives a unifying perspective of first-order relaxations based on a layer-wise convex relaxation framework, and shows that the optimal layer-wise convex relaxation exhibits a non-trivial gap. 
This gap highlights the inability of first-order relaxations to capture the ReLU activation, the complementarity of which is naturally represented with second-order (i.e., quadratic) constraints.

An approach for addressing these limitations is to lift to a second-order space where the quadratic ReLU constraints can be linearized via the reformulation linearization technique (RLT) \cite{SheraliAdams1998}. 
 \cite{RaghunathanSteinhardtEtAl2018} first proposed representing ReLUs with a quadratic constraint and lifting the problem to a SDP. 
Work in this area is primarily based on SDPs both from the primal perspective \cite{RaghunathanSteinhardtEtAl2018,DathathriDvijothamEtAl2020,AndersonMaEtAl2021,MaSojoudi2020}, and from the dual perspective \cite{FazlyabMorariEtAl2022,NewtonPapachristodoulou2021}. 
One exception is \cite{DvijothamStanforthEtAl2020}, which uses diagonally dominant matrices to relax the SDP further and represent it as a linear program (LP). 

One challenge specific to second-order relaxations is that the quantities of interest, i.e., neuron values, are obfuscated in the relaxation formulation. 
This makes it difficult to tell how the constraints interact, and to what extent they accurately capture the neural network computation.
The analysis provided in this work aims also to clarify these points.

\section{PRELIMINARIES}\label{sec:preliminaries}
\subsection{Notation and Terminology}
For a matrix $M$, we use $M_{i, j}$ to denote the entry in the $i$th row and $j$th column, $M_{i, *}$ denotes the entire $i$th row, and $M_{*, j}$ denotes the entire $j$th column.
For sets of matrices, $S_M + S_N$ denotes the Minkowski sum: $S_M + S_N : = \{ M + N \mid M \in S_M, N \in S_N\}$. 
The matrix inner product $\innerMat{A}{B}$ is defined as $\innerMat{A}{B} = \text{Tr}(A^\top B)$. 
$S^+$ is the set of positive semidefinite (PSD) matrices; $\mathcal{N}$ is the set of entrywise non-negative square matrices (we use $S^+_n$ and $\mathcal{N}_n$ to specify the dimension if it is unclear from context). 
The cone of doubly non-negative matrices is defined as $S^+ \cap \mathcal{N}$, where $\cap$ denotes the set intersection.

\subsection{Problem Setting}
\paragraph{Feedforward ReLU Networks.}
We consider an n-layer feedforward ReLU network representing a function $f$ with input $z_{0, *} \in \mathbb{R}^{h_0}$ and output $f(z_{0, *}) = z_{n, *} \in \mathbb{R}^{h_n}$, with $f$ being the composition of layer-wise functions, i.e., $f = f_{n} \circ f_{n - 1} \circ \cdots \circ f_{1}.$
The $i$th layer of $f$ corresponds to a function $f_i : \mathbb{R}^{h_{i - 1}} \rightarrow \mathbb{R}^{h_i}$ ($h_i$ is the dimension of the hidden variable $z_{i, *}$) of the form
\begin{equation}
    z_{i, *} = f_i(z_{i - 1, *}) = \sigma(\overline{W}[i - 1] z_{i - 1, *}+ b[i - 1]).
\end{equation}
where $\overline{W}[i - 1] \in \mathbb{R}^{h_i \times h_{i - 1}}$ is the weight matrix, $b[i - 1] \in \mathbb{R}^{h_i}$ is the bias.
The vector of all neurons in the $i$-th layer is denoted by $z_{i, *}$, while $z_{i, j}$ refers specifically to the $j$-th neuron in the $i$-th layer. The function $\sigma(\hat{z}_{i, j}) = \max(0, \hat{z}_{i,j}) $ is the Rectified Linear Unit (ReLU) function; on vectors, $\sigma(\cdot)$ operates element-wise.
We assume that this activation is not applied in the last layer $f_n$.
We  also model the preactivation values $\hat{z}_{i, *} = \overline{W}[i - 1] z_{i - 1, *}+ b[i - 1]$, so that $z_{i, *} = \sigma (\hat{z}_{i, *})$ for $i \geq 1$, and make the identification $\hat{z}_{0, *} = z_{0, *}$.

In this paper, we use a positive/negative splitting for each neuron, denoted by $\lambda^+\geq 0, \, \lambda^-\geq 0$, rather than the typical pre/post-activation splitting. 
The positive/negative and pre/post-activation splittings are related by the following equations:
\begin{align}
    \lambda^+_{i, j} &= z_{i, j}  \\
    \lambda^-_{i, j} &= z_{i, j} - \hat{z}_{i, j}
\end{align}

We will often need to make a number of routine conversions, chiefly between matrices that act on post-activation neurons in single layer, $z_{i, *}$, to those that act on the collection of all positive/negative splittings of variables.
To emphasize that these conversions maintain equivalence, we refer to these converted matrices using the same letter and use an overline (e.g., $\overline{M}$) to denote matrices that act on single layers, and an unmarked matrix (e.g., $M$) to denote matrices that act on the collection of all variables. 
For inequalities, this conversion also includes the addition of a slack variable to convert to an equality constraint. 
We use an underline (e.g., $\underline{M}$) to denote that a slack variable has not been included and the inequality should be treated as such. Concrete examples illustrating these conversions are included in the Appendix.

\paragraph{Safety Rules.}We consider safety rules specified as input-output relationships on $f$. Specifically, for all inputs from the set $\mathcal{X} \subseteq  \mathbb{R}^{h_0}$, we aim to ensure that the output belongs to the set $\mathcal{Y} \subseteq \mathbb{R}^{h_n}$. 
We consider bounded, polytopic input sets
\begin{equation}{\footnotesize
    \mathcal{X} = \{x \in \mathbb{R}^{h_0} \mid \overline{A} x \leq a \},}
\end{equation}

and output sets specified by a half-space constraint, 
\begin{equation}{\footnotesize
    \mathcal{Y} = \{y \in \mathbb{R}^{h_n} \mid \overline{c}^\top y \geq d \}.}
\end{equation}
This is not restrictive, as this construction can be extended to polytopic output sets as well; because polytopes are the intersection of half-spaces, each inequality defining the polytope can be verified independently.
In this paper, we consider the formulation of finding output violations given input constraints. However, the methods developed in this paper can be readily adapted to finding minimal adversarial disturbances by imposing constraints on the output and optimizing over the input.
\paragraph{Optimization Formulation.} Verification can be posed as a non-convex optimization problem:
\begin{mini}|s|
{\lambda^+,\, \lambda^-}{\overline{c}^\top(\lambda^+_{n, *} - \lambda^-_{n, *})}
{\label{opt:linear_exact}}{\texttt{OPT} =}
\addConstraint{\overline{A} (\lambda^+_{0, *} - \lambda^-_{0, *}) \leq a}
\addConstraint{\lambda^+_{i + 1, *} - \lambda^-_{i + 1, *} = \overline{W}[i]\lambda^+_{i, *} + b[i]} 
\addConstraint{\lambda^+_{i, j}\lambda^-_{i, j} = 0 \quad \forall i, j}
\addConstraint{\lambda^+ \geq 0,\, \lambda^- \geq 0}
\end{mini}
where the network is safe if and only if $\texttt{OPT} \geq d$.

\subsection{Completely Positive Programming}
Completely positive programs (CPP) are linear optimization problems over matrix variables \cite{Dure2010}:
\begin{mini}
{X}{\innerMat{M_0}{X}}
{\label{opt:generic_CPP}}{}
\addConstraint{\innerMat{M_i}{X} = m_i, \quad i \in \{1, \dots, L\}}
\addConstraint{X \in C^*}
\end{mini}
where $C^* \subset S^+ \cap \mathcal{N}$ is the cone of completely positive (CP) matrices, that is, matrices that have a factorization with entrywise non-negative entries:
\begin{equation}
{\footnotesize
    \mathcal{C}^*_{n} := \{X \in \mathbb{R}^{n \times n} \mid X = \sum_k x^{(k)} (x^{(k)})^\top, \quad x^{(k)} \in  \mathbb{R}^n_{\scriptscriptstyle \geq 0} \}}
\end{equation}
Based on the Sum of Squares (SOS) hierarchy, \cite{Parrilo2000} constructed a hierarchy of cones $\{(\mathcal{K}^r)^*\} $ approximating the completely positive cone from the exterior. Meaning, $\mathcal{C}^* = \bigcap_{r\geq 0} (\mathcal{K}^r)^*$, and
\begin{equation}\label{eq:SOS_hierarchy}{\footnotesize
S^+ \cap \mathcal{N} = (\mathcal{K}^0)^* \supset (\mathcal{K}^1)^* \supset \ldots }
\end{equation}
Optimizing over each cone, $(\mathcal{K}^r)^*$, can be posed as an SDP. In addition, the hierarchy gives a formulaic methodology for constructing relaxations of CPPs that trade off accuracy and tractability.

\section{GEOMETRY OF COMPLETE POSITIVITY} \label{sec:geometry}
Before presenting the full CPP formulation, this section will provide an introductory exploration of the interaction between complete positivity and the verification problem. 
We focus on a single neuron and provide intuition for the geometry of CP matrices when representing a ReLU activation.\footnote{Although, for dimensions $n \leq 4$, $C^*_n = S^+_n \cap \mathcal{N}_n$.}
A single neuron is represented by two variables, $\lambda^+ \geq 0$ and $\lambda^- \geq 0$, denoting the positive and negative split of the neuron, respectively. 
The ReLU function is then defined by the quadratic constraint $\lambda^+ \lambda^- = 0$.
The completely positive formulation entails lifting to a matrix, 
\begin{equation}\label{eq:single_neuron_cp}
    {\scriptsize
    \begin{pmatrix}
    \Lambda[\lambda^+, \lambda^+]  & \Lambda[\lambda^+, \lambda^-]  &\lambda^+ \\
    \Lambda[\lambda^+, \lambda^-]  & \Lambda[\lambda^-, \lambda^-]  & \lambda^- \\
    \lambda^+ & \lambda^- & 1 
    \end{pmatrix}
    \in C^*}
\end{equation} 
where the elements of $\Lambda \in \mathbb{R}^{2 \times 2}$ correspond to cross terms of $\lambda^+$ and $\lambda^-$, and are denoted with symbolic indexing.
In the lifted formulation, the ReLU is further defined by the linear constraint
$
    \Lambda[\lambda^+, \lambda^-] = 0.
$
Complete positivity means there exists a factorization
\begin{equation}\label{eq:single_factor} {\scriptsize
    \begin{pmatrix}
        \Lambda[\lambda^+, \lambda^+]  & \Lambda[\lambda^+, \lambda^-]  &\lambda^+ \\
        \Lambda[\lambda^+, \lambda^-]  & \Lambda[\lambda^-, \lambda^-]  & \lambda^- \\
        \lambda^+ & \lambda^- & 1 
    \end{pmatrix}
    =
    \sum_{k=1}^K
    \begin{pmatrix} \lambda^{+, (k)} \\ \lambda^{-, (k)} \\ \xi^{(k)} \end{pmatrix}
    \begin{pmatrix} \lambda^{+, (k)} \\ \lambda^{-, (k)} \\ \xi^{(k)} \end{pmatrix}^\top}
\end{equation}
where $\lambda^{+, (k)}, \, \lambda^{-, (k)}, \,  \xi^{(k)} \geq 0$, and $K$ is the rank of the factorization.
The variables in CPPs and SDPs are weighted (convex in the case of CPPs) combinations of the factors, e.g., neurons $\lambda^\pm = \sum_k \xi^{(k)} \lambda^{\pm, (k)}$.
As we will further explore in Section \ref{sec:results}, it is important that each factor satisfies the constraints, however we are only able to impose them on their aggregate.
A relaxation gap arises when constraints that hold for the sum of factors do not hold for individual factors, e.g., $\Lambda[\lambda^+, \lambda^-] = \sum_k \lambda^{+, (k)} \lambda^{-, (k)} = 0$, but  $\lambda^{+, (k)} \lambda^{-, (k)} \neq~0$.

\begin{figure}
    \centering
    \includegraphics[width=0.75\columnwidth]{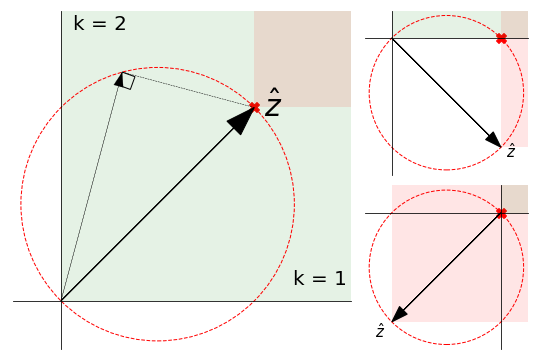}
    \caption{This figure depicts the geometry of CP matrices and ReLU constraints. $\sum_k z^{(k)} (\hat{z}^{(k)} - z^{(k)}) = 0$ constrains $z$ to the dotted red circle,
    $z^{(k)} \geq \hat{z}^{(k)}$ is depicted by the shaded red region, and $z^{(k)} \geq 0$ is depicted by the shaded green region. $z$ (indicated by the red cross) is fully determined by these constraints.}
    \label{fig:cpp_geometry}
\end{figure}

In the context of ReLUs, non-negativity of $\lambda^{+, (k)}, \, \lambda^{-, (k)}, \,  \xi^{(k)}$ is critical for ensuring that each factor corresponds to a valid ReLU, even if $K \neq 1$.
In terms of the pre/post-activation splitting,
\begin{equation}\label{eq:cp_factor_constraints}{\footnotesize
    z^{(k)} \geq 0, \quad z^{(k)} \geq \hat{z}^{(k)}, \quad \sum_{k=1}^K z^{(k)} (z^{(k)} - \hat{z}) = 0.}
\end{equation}
Figure \ref{fig:cpp_geometry} depicts three different cases when $K = 2$. In the diagram, each factor satisfies the ReLU constraints, i.e., $z^{(k)} = \sigma(\hat{z}^{(k)})$.  Intuitively, the exactness of the CPP formulation follows from analogous arguments that each term in the factorization of an optimal solution is feasible for \eqref{opt:linear_exact}; this notion is formalized in the proof of Theorem~\ref{thm:verif_cpp}.

\section{RESULTS}\label{sec:results}
In this section, we first state our main result, the exact reformulation of \eqref{opt:linear_exact} as a CPP, in Theorem \ref{thm:verif_cpp}. In order to understand SDP relaxations of this reformulation, the remainder of this section is dedicated to breaking down the constraints in the CPP and illustrating how the guarantees change when the CP constraint is relaxed. The proofs of all results stated in this section are provided in the Supplemental Material.

In this section, the objects of interest will be matrices of the form $\quadmat{\Lambda}{\lambda}{1}$, where $\Lambda$ is a square matrix representing cross terms of elements in $\lambda$, a vector. 
In the proposed formulation, $\lambda^\top = \begin{pmatrix} (\lambda^+)^\top & (\lambda^-)^\top & s^\top \end{pmatrix}$ is the concatenation of $\lambda^+, \, \lambda^-$, the positive and negative splitting for each neuron, and $s$, slack variables used to represent inequality constraints.
To emphasize the relationship with $\lambda$, we use symbolic indexing to refer to specific terms in $\Lambda$; for example, $\Lambda[\lambda^+_{i, j}, s_k]$ refers to the element in $\Lambda$ corresponding to the dot product of factors corresponding to $\lambda^+_{i, j}$ and $s_k$. 

\begin{theorem}\label{thm:verif_cpp}
    Problem \eqref{opt:linear_exact} is equivalent to the following convex optimization problem:
    \begin{mini!}|s|
    {\lambda,\, \Lambda}{c^\top \lambda}
    {\label{opt:CPP}}{\textup{\texttt{OPT}}_{\textup{\texttt{CPP}}} =}
    \addConstraint{A_{j, *} \lambda = a_j  \label{eq:input_lin}}
    \addConstraint{\innerMat{A_{j, *}^\top A_{j, *}}{\Lambda} = a^2_j \label{eq:input_self_quad}}
    \addConstraint{W[i]_{j, *} \lambda = b[i]_j \label{eq:nw_lin}}
    \addConstraint{\innerMat{W[i]_{j, *}^\top W[i]_{j, *}}{\Lambda} = b[i]^2_j \label{eq:nw_self_quad}}
    \addConstraint{\Lambda[\lambda^+_{i, j}, \lambda^-_{i, j}] = 0 \label{eq:relu}}
    \addConstraint{\quadmat{\Lambda}{\lambda}{1} \in \mathcal{C}^* \label{eq:comp_pos}}
    \end{mini!}
    
    Exactness is defined by the following conditions: 
    \begin{enumerate}
        \item $\textup{\texttt{OPT}} =  \textup{\texttt{OPT}}_{\textup{\texttt{CPP}}}$: Problem \eqref{opt:linear_exact} and Problem \eqref{opt:CPP} have the same objective value.
        \item If $(\lambda^*, \Lambda^*) = \arg \min \eqref{opt:CPP}$ then $\lambda^*$ is in the convex hull of optimal solutions for Problem \eqref{opt:linear_exact}.
    \end{enumerate}
\end{theorem}

    The salient features of \eqref{opt:CPP} are the \emph{pairs} of linear and ``self-quadratic" constraints and the CP constraint.

    To motivate the following discussion, we introduce the zeroth order sum of squares (0-SOS) relaxation as a natural SDP relaxation for Problem \eqref{opt:CPP}:
    \begin{mini}
    {\lambda, \, \Lambda}{c^\top \lambda}
    {\label{opt:0SOS}}{}
    \addConstraint{\eqref{eq:input_lin}-\eqref{eq:relu}}
    \addConstraint{\quadmat{\Lambda}{\lambda}{1} \in S^+ \cap \mathcal{N}}
    \end{mini}
    
    This relaxation (derived from \eqref{eq:SOS_hierarchy}) replaces the intractable CP constraint, $\quadmat{\Lambda}{\lambda}{1} \in C^*$ , with a tractable doubly non-negative constraint, $\quadmat{\Lambda}{\lambda}{1} \in S^+ \cap \mathcal{N}$ . 
    
    Next, we explore the cooperative role that the linear and self-quadratic constraints play in enforcing linear constraints on individual factors. In Theorem \ref{thm:lin_eq_factor}, we highlight the identical role they play in Problems \eqref{opt:CPP} and \eqref{opt:0SOS}  when representing equality constraints. They do, however, result in different guarantees when representing inequality constraints. These differences are explored in Corollary \ref{cor:ineq_cpp_0sos}.
    We present results in a generic form, because none of these results rely on the specific form of the verification problem, and it allows us to focus on the core assumptions needed to derive each result. 
    To tie things back to \eqref{opt:CPP} and \eqref{opt:0SOS}, we intersperse these results with discussion of implications for verification.

    The following theorem shows that pairs of linear and self-quadratic constraints play an identical role in both the CPP and SDP-based formulations.
    \begin{theorem}\label{thm:lin_eq_factor}
            Suppose there is a factorization
    \begin{equation}{\small
        \quadmat{\Lambda}{\lambda}{1} 
        = \sum_{k = 1}^K 
        \begin{pmatrix} \lambda^{(k)}\\ \xi^{(k)} \end{pmatrix}
        \begin{pmatrix} \lambda^{(k)}\\ \xi^{(k)} \end{pmatrix}^\top}.
    \end{equation}
    Then, $ V_{i, *} \lambda^{(k)} = \xi^{(k)} v_i$ for all $k$ if and only if
    \begin{align}\footnotesize
        V_{i, *}\lambda &= v_i \label{eq:lin_generic} \\
        \innerMat{V_{i, *}^\top V_{i, *}}{\Lambda} &= v_i^2. \label{eq:self_quad_generic}
    \end{align}
    \end{theorem}
    This theorem underscores the importance of \emph{pairing} linear and self-quadratic constraints to enforce linear constraints on the factors, $\lambda^k$ and $\xi^k$. Notice that this result does not rely on complete positivity, and holds even with a weaker positive semidefinite constraint.
    
    \paragraph{Network-Defining Constraints.} Theorem \ref{thm:lin_eq_factor} characterizes how constraints \eqref{eq:nw_lin} and \eqref{eq:nw_self_quad} act in Problems \eqref{opt:CPP} and \eqref{opt:0SOS}. They collectively ensure that the neurons in consecutive layers are related through the weights and biases of the networks.

   Because the equality constraints act the same way in \eqref{opt:CPP} and \eqref{opt:0SOS}, the only potential source of a relaxation gap is due to inequality constraints. 
   In our framework, we convert inequality constraints to equality constraints through the introduction of slack variables. 
   For a CP matrix, each factor's slack variables are guaranteed to be non-negative, whereas for doubly non-negative matrices, only their weighted sum is guaranteed to be non-negative. 
   That is, individual slack variables may be negative, and thus individual factors are not guaranteed to satisfy the inequalities. The following corollary highlights these differences. 
     
\begin{corollary} \label{cor:ineq_cpp_0sos}
    Suppose there is a factorization,
    \begin{equation} \label{eq:slack_factor}{\scriptsize
        \begin{pmatrix}
        \Lambda[\lambda, \lambda^\top] & \Lambda[\lambda, s^\top] & \lambda\\
        \Lambda[s, \lambda ^\top] & \Lambda[s, s^\top] & s\\
        \lambda^\top & s^\top & 1
        \end{pmatrix}
         = \sum_{k = 1}^K 
        \begin{pmatrix} \lambda^{(k)}\\ s^{(k)} \\\xi^{(k)} \end{pmatrix}
        \begin{pmatrix} \lambda^{(k)}\\ s^{(k)} \\\xi^{(k)} \end{pmatrix}^\top}
    \end{equation}
    
    and for all $i = \{1, \ldots L\}$, $(\lambda, \, \Lambda)$ satisfy
    \begin{align} \footnotesize
        &V_{i, *}\lambda + s_i = v_i, \label{eq:slack_lin}\\
        &\innerMat{V^\top_{i, *} V_{i, *}}{\Lambda[\lambda, \lambda^\top]} + 2 V_{i, *} \Lambda[\lambda, s_i] + \Lambda[s_i, s_i] = v_i^2. \label{eq:slack_quad}
    \end{align}

    Then we must have:
    \begin{enumerate}
        \item For $\scriptsize \begin{pmatrix}
        \Lambda[\lambda, \lambda^\top] & \Lambda[\lambda, s^\top] & \lambda\\
        \Lambda[s, \lambda ^\top] & \Lambda[s, s^\top] & s\\
        \lambda^\top & s^\top & 1
        \end{pmatrix} \in C^*$,
        each factor individually satisfies the inequalities:
        \begin{equation}{\footnotesize
            V_{i, *}\lambda^{(k)} \leq \xi^{(k)} v_i.}
        \end{equation}
        \item For $\scriptsize  \begin{pmatrix}
        \Lambda[\lambda, \lambda^\top] & \Lambda[\lambda, s^\top] & \lambda\\
        \Lambda[s, \lambda ^\top] & \Lambda[s, s^\top] & s\\
        \lambda^\top & s^\top & 1
        \end{pmatrix} \in S^+$, and $s \geq 0$
        the weighted sum of factors satisfies the inequalities:
        \begin{equation}{\footnotesize
            V_{i, *} \lambda = \sum_k \xi^{(k)} V_{i, *}  \lambda^{(k)} \leq v_i}.
        \end{equation}
    \end{enumerate}
\end{corollary}
    \paragraph{Input and Non-negativity Constraints.} This corollary shows that constraints \eqref{eq:input_lin} and \eqref{eq:input_self_quad} are sufficient for ensuring that each factor satisfies the input constraints in \eqref{opt:CPP}. However, we cannot derive the same result for \eqref{opt:0SOS} because membership in the cone $S^+ \cap \mathcal{N}$ does not necessarily imply the existence of a non-negative factorization. An analogous conclusion applies to the non-negativity of $\lambda^+$ and $\lambda^-$. 
    
    \paragraph{ReLU Constraints.} In \eqref{opt:CPP}, each factor satisfies $\lambda_{i, j}^{\pm, (k)}\geq 0$. Thus the constraint $\sum_{k} \lambda_{i, j}^{+, (k)}\lambda_{i, j}^{-, (k)} = 0 $ can only hold if $\lambda_{i, j}^{+, (k)}\lambda_{i, j}^{-, (k)} = 0$ for all $k$, i.e., each factor satisfies the ReLU constraints. In \eqref{opt:0SOS}, it is only guaranteed that $\sum_k \xi^{(k)} \lambda_{i, j}^{\pm, (k)}\geq 0$, so the same conclusion cannot be drawn.

    \paragraph{Verification-Defining Constraints.} In the proposed CPP formulation, constraints \eqref{eq:input_lin}-\eqref{eq:relu} are critical for encoding the verification problem. By the results presented in this section, the removal of any of these constraints would fundamentally misrepresent the problem. For this reason, we call them the verification-defining constraints.
    
    \paragraph{Strengthening Constraints.}  In the (0-SOS) relaxation \eqref{opt:0SOS}, however, the verification-defining constraints are not sufficient to ensure exactness. One common strategy is to impose additional ``strengthening constraints." These are constraints encoding \emph{a priori} knowledge of the optimal solution. They will have no effect on \eqref{opt:CPP}, but can potentially reduce the feasible domain in relaxations, thereby generating tighter relaxations. In the context of verification, these are typically derived from bound constraints:

    \noindent\begin{tabularx}{\linewidth}{@{}X@{\hspace{-1em}}X@{}}
        \vspace{-1cm}
        \begin{equation}{\footnotesize\label{eq:hat_ub}
            \!\!\!\hat{u}_{i, j} - (\lambda^+_{i, j} - \lambda^-_{i, j}) \geq 0}
        \end{equation}
        \vspace{-0.75cm}
        &
        \vspace{-1cm}
        \begin{equation}{\footnotesize\label{eq:plain_ub}
            u_{i, j} - \lambda^+_{i, j} \geq 0 }
        \end{equation}
        \vspace{-0.75cm} \\
        \vspace{-0.75cm}
        \begin{equation}{\footnotesize\label{eq:hat_lb}
            (\lambda^+_{i, j} - \lambda^-_{i, j}) - \hat{l}_{i, j} \geq 0}
        \end{equation}
        \vspace{-0.25cm} &
        \vspace{-0.75cm}
        \begin{equation}{\footnotesize\label{eq:plain_lb}
            \lambda^+_{i, j} - l_{i, j} \geq 0}
        \end{equation}
        \vspace{-0.25cm}
    \end{tabularx}
    These are lower and upper bounds for each neuron's pre/post-activation values (denoted by $\hat{l}_{i, j}, \hat{u}_{i, j}$, and $l_{i, j}, u_{i, j}$ respectively)---they can be found efficiently by forward propagating bounds on the input set. These bounds can be further refined by leveraging the following characterization of the convex hull of the graph of ReLUs:
    \begin{equation}\label{eq:triangle_cut}{\footnotesize
        \lambda^+_{i, j} \leq \frac{u_{i, j} - l_{i, j}}{\hat{u}_{i, j} - \hat{l}_{i, j}} ((\lambda^+_{i, j} - \lambda^-_{i, j}) - \hat{l}_{i, j} )} + l_{i, j}
    \end{equation}
    This inequality is commonly referred to as the triangle relaxation due to the geometry of the convex hull \cite[Fig. 6.4]{LiuArnonEtAl2021}.
    
    In our framework, each of these inequality constraints would be imposed by introducing a slack variable and including a pair of linear and self-quadratic constraints. In principle, \eqref{eq:hat_ub}-\eqref{eq:triangle_cut} could be further combined with additional valid inequalities to provide even stronger relaxations.
    However, the following Corollary shows that simply combining existing linear inequalities, e.g., by multiplication, will not actually strengthen the relaxation.
    \begin{corollary}\label{cor:cross_redundant}
        If
        $\scriptsize \begin{pmatrix}
        \Lambda[\lambda, \lambda^\top] & \Lambda[\lambda, s^\top] & \lambda\\
        \Lambda[s, \lambda ^\top] & \Lambda[s, s^\top] & s\\
        \lambda^\top & s^\top & 1
        \end{pmatrix} \in S^+ \cap \mathcal{N}$, 
            and \eqref{eq:slack_factor}, \eqref{eq:slack_lin}, \eqref{eq:slack_quad} hold,
            then the following holds:   
        \begin{equation}\label{eq:cross_quad_constr}{\footnotesize 
            v_i v_j - v_i V_{j, *} \lambda  - v_j V_{i, *}\lambda  + \innerMat{V^\top_{i, *} V_{j, *}}{\Lambda[\lambda, \lambda^\top]} \geq 0}
        \end{equation}
    \end{corollary}
    Inequality \eqref{eq:cross_quad_constr} is the linearized inequality representing $(v_i - V_{i, *} \lambda)(v_j - V_{j,*}\lambda) \geq 0$. In contrast to the self-quadratic constraints, this is a product of two different linear constraints. For this reason, we call constraints of this form ``cross-quadratic" constraints.
    Corollary \ref{cor:cross_redundant} states that if we enforce constraints of the form $V_{i, *}\lambda \leq v_i$ through the introduction of a slack variable and paired linear/self-quadratic constraints, then the derived cross-quadratic constraint will also hold. 
    This is significant because in a number of existing works, cross-quadratic constraints derived from \eqref{eq:hat_ub}-\eqref{eq:plain_lb} are introduced, however the linear and self-quadratic constraints are excluded.

\section{UNIFYING SDP-BASED VERIFICATION} \label{sec:existing}
  
  Now that we have established the notions of verification-defining and strengthening constraints, we are in a position to show how existing work fits into our framework. 
  We discuss two seemingly disparate frameworks that have served as the basis for subsequent work, \cite{RaghunathanSteinhardtEtAl2018} and \cite{FazlyabMorariEtAl2022}. 
  In this section we show how these frameworks, as well as subsequent works, are related by couching their approaches in the language of our proposed verification-defining and strengthening constraints. 
  A number of subsequent works focus on improving computational aspects of existing relaxations (e.g., time to solution, memory usage) without aiming to improve the relaxation gap.
  These works are discussed in the Appendix rather than the main body of this paper because they largely have the same constraint structure at work.
  
  A number of works in this section construct strengthening cuts derived from pre/post-activation bounds on neurons. 
  Without loss of generality, we will assume that $\hat{l}_{i, j} < 0 < \hat{u}_{i, j}$ (otherwise the ReLU could be treated as the identity, if $\hat{l}_{i, j} \geq 0$, or zero, if $\hat{u}_{i, j}\leq 0$). The post-activation bounds are then $l_{i, j} = 0$ and $u_{i, j} = \hat{u}_{i, j}$.

\begin{figure}
    \centering
    \includegraphics[width=1.0\linewidth]{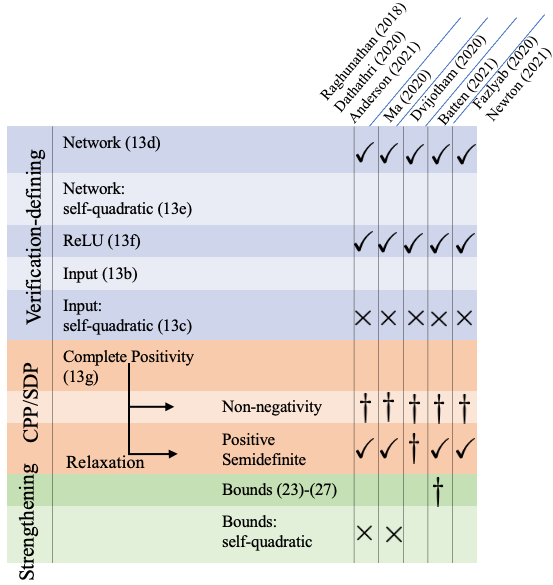}
    \caption{This table summarizes the constraints included in existing work. The $\dagger$ indicates partial encoding of constraints/categorizations that require additional qualification. Cross-quadratic constraints are marked with a $\times$ in the self-quadratic row.}
    \label{fig:chart}
\end{figure}

\subsection{Direct SDP}
    To the best of our knowledge, \cite{RaghunathanSteinhardtEtAl2018} was the first to propose an SDP relaxation of the verification problem:
    \begin{mini}|s|
    {z, Z}{\overline{c}^\top z_{n, *}}
    {}{\label{opt:raghunation_raw}}
    \addConstraint{z_{i, *} \geq 0, \, z_{i, *} \geq \overline{W}[i - 1] z_{i - 1, *}}
    \addConstraint{Z[z_{i, j}, z_{i, j}] = \overline{W}[i - 1]_{j, *} Z[z_{i - 1, *}, z_{i,j}]}
    \addConstraint{-Z[z_{i,j}, z_{i,j}] + (u_{i, j} + l_{i, j})z_{i,j} - l_{i, j}u_{i,j} \geq 0, \, \forall i, j}
    \addConstraint{\quadmat{Z}{z}{1} \in \mathcal{S}^+}
    \end{mini}
    Although their approach does not include biases, it is straightforward to include them. Problem \eqref{opt:raghunation_raw} is equivalent to the following optimization problem:
    \begin{mini}|s|
    {\lambda, \, \Lambda}{c^\top \lambda}
    {}{\label{opt:SDR}}
    \addConstraint{(l_{0, i} + u_{0, i}) U_{i, *}\lambda - \innerMat{U^\top_{i, *}U_{i, *}}{\Lambda} - u_il_i \geq 0}
    \addConstraint{\hat{u}_{i, j} \lambda^+_{i, j} -  \Lambda[\lambda^+_{i, j}, \lambda^+_{i, j}] \geq 0, \, i \geq 1}
    \addConstraint{W[i]_{j, *} \lambda = b[i]_j}
    \addConstraint{\Lambda[\lambda^+_{i, j}, \lambda^-_{i, j}] = 0}
    \addConstraint{\quadmat{\Lambda}{\lambda}{1} \in S^+}
    \addConstraint{\lambda \geq 0}
    \end{mini}
    where $U$ is defined as $U_{i, *} \lambda = \lambda^+_{0, i} - \lambda^-_{0, i}$.
    In the context of \eqref{opt:0SOS}, we see that \eqref{opt:SDR} does not account for non-negativity of $\Lambda$, i.e., $\Lambda \in \mathcal{N}$, and the network self-quadratic constraints \eqref{eq:nw_self_quad}. \eqref{opt:SDR} also includes strengthening constraints derived from  $\lambda^+_{i, j} \geq 0$ and $(\hat{u}_{i, j} - \lambda^+_{i, j}) \geq 0$.

    \paragraph{Tightening Extensions.}
        Based on the observation that \eqref{opt:SDR} sometimes produces bounds looser than the triangle LP relaxation , \cite{BattenKouvarosEtAl2021} proposes including linear cuts derived from \eqref{eq:triangle_cut}. In the following, we provide intuition for the gap that the triangle inequality closes.
        
        \begin{figure}
            \centering
            \includegraphics[width=\columnwidth]{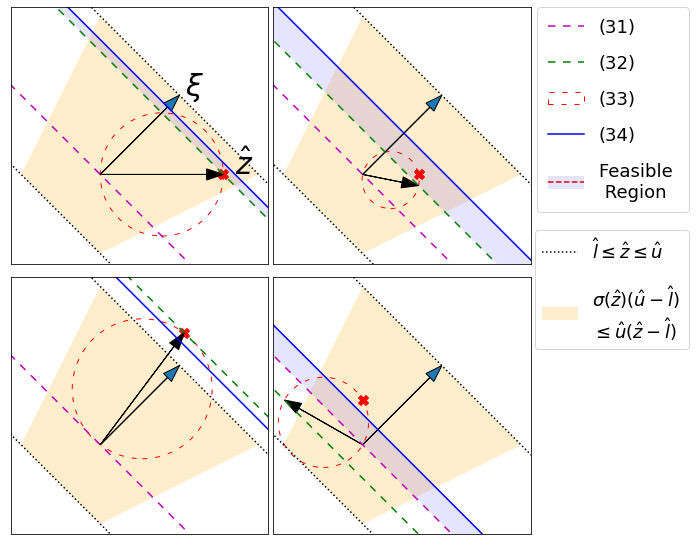}
            \caption{This figure depicts the geometry of the triangle cut. $\sum_k z^{(k)} (\hat{z}^{(k)} - z^{(k)}) = 0$ constrains $z$ to the dotted red circle. $z \geq \hat{z}$ and $z \geq 0$ constrain $z$ to the right of the dotted green and purple lines, respectively. The triangle cut~\eqref{eq:tri_ub} constrains z to the left of the solid blue line. The lower and upper bounds ($\hat{l} \leq \hat{z}$ and $\hat{z} \leq \hat{u}$) used to construct the triangle cut are demarcated by the dotted black lines. The red cross marks the factor-wise ReLU (as in Figure~\ref{fig:cpp_geometry}), and the shaded orange region indicates the region (in terms of $\hat{z}$) for which the factor-wise ReLU is feasible.}
            \label{fig:triangle_cut_geometry}
        \end{figure}
        
        The PSD constraint ensures that each neuron has a factorization in the form of \eqref{eq:single_factor} where $\lambda^{+, (k)}, \, \lambda^{-, (k)}$ are not necessarily non-negative (without loss of generality $\xi^{(k)}$ can be normalized to be non-negative). In contrast to the CP case, the constraint $\lambda^+ \geq 0$ only guarantees that $\sum_{i = 1}^K \xi^{(k)}\lambda^{+, (k)} \geq 0$, and similarly for other inequalities. In terms of the pre/post-activation splittings, it is guaranteed that:
        \begin{align}
        &\label{eq:inact_lb}\medmath{ \sum_{i = 1}^K \xi^{(k)}z^{(k)} \geq 0},\\
        &\label{eq:act_lb}\medmath{ \sum_{i = 1}^K \xi^{(k)}z^{(k)}\geq \sum_{i = 1}^K \xi^{(k)}\hat{z}^{(k)}},\\
        & \label{eq:compl}\medmath{\sum_{i = 1}^K z^{(k)}(\hat{z}^{(k)} - z^{(k)})= 0.}
        \end{align}
        These are the analogs of \eqref{eq:cp_factor_constraints}.
        The triangle cut provides the additional guarantee that
        \begin{equation}
        \label{eq:tri_ub} \medmath{\sum_{i = 1}^K \xi^{(k)}z^{(k)} \leq \frac{\hat{u}_{i, j} }{\hat{u}_{i, j} - \hat{l}_{i, j}} \left(\sum_{i = 1}^K \xi^{(k)}\hat{z}^{(k)}  - \hat{l}_{i, j} \right).}
        \end{equation}
        Figure \eqref{eq:triangle_cut} illustrates the geometry of these constraints when $K = 2$.
        The plots, from left to right and top to bottom, correspond to the cases where \eqref{eq:tri_ub} is (i) non-redundant and includes $\sigma(\hat{z})$, (ii) redundant, (iii) infeasible, and (iv) non-redundant but eliminates $\sigma(\hat{z})$ given \eqref{eq:inact_lb} - \eqref{eq:compl}.
        The diagrams in the top row illustrate the effect of the triangle cut in restricting the feasible region of $z$ given $\hat{z}$.
        In the first case, the triangle cut strictly reduces the feasible region, resulting in a strengthened relaxation. 
        In the second case, the region defined by the intersection of \eqref{eq:inact_lb} - \eqref{eq:compl} lies entirely within the region defined by \eqref{eq:tri_ub}, i.e., the triangle cut does not add any additional information.
        
        The diagrams in the bottom row illustrate the triangle cut's role in enforcing ``acceptable" values of $\hat{z}$.
        The third (infeasible) case occurs if and only if the pre-activation bounds are actually not satisfied $\sum_{i = 1}^K \xi^{(k)}\hat{z}^{(k)} \not \in [\hat{l}, \hat{u}]$, i.e., the triangle cut implicitly constrains $\hat{z} \in [\hat{l}, \hat{u}]$.
        In the fourth case, the triangle cut is feasible but counter-intuitively renders the factor-wise ReLU, $z^{(k)} = \sigma(\hat{z}^{(k)})$, infeasible.
        This is because the triangle cut is only valid for $\hat{z}$ that are convex combinations of  
        \begin{equation}\label{eq:tri_cond}
        \small{\hat{z}^{(k)} \in \left[\frac{\hat{l}}{\xi^{(k)}}, \frac{\hat{u}}{\xi^{(k)}}\right],}
        \end{equation}
        even if $\sum_{i = 1}^K \xi^{(k)}\hat{z}^{(k)} \in [\hat{l}, \hat{u}]$. 
        Condition \eqref{eq:tri_cond} is satisfied by any $z^{(k)}$ corresponding to a forward pass of the network, but cannot be explicitly enforced in an SDP.
        The behavior of the triangle cut illustrates the contrast between the guarantees of the CPP formulation, which can be leveraged to (implicitly) impose constraints on individual factors, and SDP relaxations, which can only enforce constraints on the aggregate of factors.

    \subsection{Quadratic Constraints and the S-procedure}
    \cite{FazlyabMorariEtAl2022} proposes a novel framework based on the notion of quadratic constraints.\footnotemark\ 
    In contrast with existing work, the proposed framework deals with classes of constraints obtained by taking infinite combinations of ``atomic" constraints. 
    For example, the atomic constraints $x\geq 0, \, y\geq 0$ are subsumed by the infinite class of constraints $\{\alpha x + \beta y \geq 0\}$ parametrized by $\alpha \geq 0, \, \beta \geq 0$. 
    Critically, $\alpha x + \beta y \geq 0$ for all $\alpha \geq 0, \, \beta \geq 0$ if and only if $x \geq 0$ and $y \geq 0$, so these constructions are equivalent.
    Similarly, for each of the infinite classes of quadratic constraints proposed, we will identify the equivalent set of atomic constraints, and rewrite the formulation of \cite{FazlyabMorariEtAl2022} in a form that parallels the (0-SOS) relaxation \eqref{opt:0SOS}.

    \footnotetext{The techniques of \cite{FazlyabMorariEtAl2022} can be applied to activations other that ReLUs by abstracting non-ReLU activations with linear bounds characterizing various properties of activation functions (e.g., monotonicity, bounded slope, bounded values). In this section, we will focus on results pertaining to ReLU networks with polytopic input constraints, and halfspace output constraints, in order to draw a direct comparison to the other methods discussed in this paper. However, the analysis of linear and self-quadratic constraints applies to SDP encodings of the properties arising from general activation functions, even in the absence of exactness. }
    
    \paragraph{Polytopic Input Sets.}
    \cite{FazlyabMorariEtAl2022} proposes encoding polytopic input sets, $\mathcal{X} =\{x \in \mathbb{R}^{k_0} \mid \overline{A} x \leq a\}$, via a quadratic constraint
    \begin{equation}
        \sum_{i, j} \Gamma_{i, j} (\underline{A}_{i, *} \lambda  - a_i)(\underline{A}_{j, *}\lambda - a_j) \geq 0
    \end{equation}
    with parameters, $\Gamma_{i, j} \geq 0, \, \Gamma_{i, i} = 0$. The equivalent atomic constraints are
    \begin{equation}
        (\underline{A}_{i, *} \lambda  - a_i)(\underline{A}_{j, *}\lambda - a_j) \geq 0, \quad \forall i \neq j.
    \end{equation}

    \paragraph{ReLU Constraints.}
    The ReLU constraints are defined through a quadratic constraint
    \begin{equation}
        \begin{aligned}
        {\footnotesize 0 \geq}  &{\footnotesize \sum_{(i, j)} \left[\rho_{(i, j)} \lambda^+_{i, j}\lambda^-_{i, j} - \nu_{i, j}  \lambda^+_{i, j}  - \eta_{i, j} \lambda^-_{i, j} \right]}\\
        &{\footnotesize + \sum_{(i, j) \neq (k, l)} \gamma_{(i, j), (k, l)}(\lambda^+_{i, j} - \lambda^+_{k, l})(\lambda^-_{i, j} - \lambda^-_{k, l})}
        \end{aligned}
    \end{equation}
     for parameters $\rho_{(i, j)} \in \mathbb{R}, \, \gamma_{(i, j), (k, l)}, \, \nu_{i, j} , \, \eta_{i, j} \geq 0$. The equivalent atomic constraints are:
    \begin{align} \footnotesize
        &\lambda^+_{i, j} \geq 0, \quad &&\lambda^-_{i, j} \geq 0, \quad \lambda^+_{i, j} \lambda^-_{i, j} = 0 \quad  \forall (i, j)\label{eq:qc_relu}\\
        &\lambda^+_{i, j}\lambda^-_{k, l} \geq 0, \quad &&\lambda^-_{i, j} \lambda^+_{k, l} \geq 0\quad \forall (i, j) \neq (k, l)\label{eq:qc_cross}
    \end{align}
    Notice that \eqref{eq:qc_relu} is equivalent to $\lambda \geq 0$ and \eqref{eq:relu}, and \eqref{eq:qc_cross} is equivalent to the non-negativity constraint $\Lambda[\lambda^+, \lambda^-] \geq 0$. 
    In summary, the formulation proposed by \cite{FazlyabMorariEtAl2022} can be shown to be equivalent to 
    \begin{mini}|s|
    {\lambda,\, \Lambda}{c^\top \lambda}
    {\label{opt:QC}}{}
    \addConstraint{
    \innerMat{\underline{A}_{i, *}^\top \underline{A}_{j, *}}{\Lambda} - a_i \underline{A}_{j, *} \lambda - a_j \underline{A}_{i, *} \lambda\geq a_ia_j, \, i \neq j}
    \addConstraint{W[i]_{j, *} \lambda = b[i]_j }
    \addConstraint{\Lambda[\lambda^+_{i, j}, \lambda^-_{i, j}] = 0 }
    \addConstraint{\quadmat{\Lambda}{\lambda}{1} \in S^+ }
    \addConstraint{\lambda \geq 0, \quad \Lambda[\lambda^+, (\lambda^-)^\top ] \geq 0}
    \end{mini}
    
    In the context of \eqref{opt:0SOS}, we see that \eqref{opt:QC} does not account for \eqref{eq:input_lin}, \eqref{eq:input_self_quad}, \eqref{eq:nw_self_quad}, and only partially accounts for $\Lambda \in \mathcal{N}$. We also observe that \eqref{opt:QC} has included cross-quadratic input constraints that are redundant in \eqref{opt:0SOS}.

    Using the 0-SOS relaxation as a shared framework, we can thus unify many of the existing works on SDP-based verification.
    Our analysis clearly lays out the constraints that are included in or excluded from existing work, with Figure \ref{fig:chart} summarizing these points.
\section{EXPERIMENTS}\label{sec:experiments}
In Sections \ref{sec:results} and \ref{sec:existing}, we have claimed that the ``verification-defining'' constraints are critical for accurately representing the neural network computation, and shown that a number of existing SDP-based relaxations omit various such constraints. In this section, we empirically evaluate some of the questions we have left unanswered. Chiefly, to what degree do the verification-defining constraints matter, i.e., is it possible to get a reasonable verification gap even without all of the verification-defining constraints? And, do proposed strengthening constraints sufficiently make up for the missing verification-defining constraints?

In this section, we consider a ReLU network with two inputs, one output, and one hidden layer with 10 neurons. We generated 100 network instances by selecting the weights and biases uniformly between -1.0 and 1.0. We consider a safety rule where the input set is defined as $\mathcal{X} := \{z_{0, *} \mid  z_{0, *} \in [-1.0, \, 0.1]^2 \}$ and the output set is defined as $\mathcal{Y}:= \{z_{n, *} \mid z_{n, *} \geq 0\}$. 

We computed the ground-truth minimum by formulating the verification problem as a MILP \cite{TjengXiaoEtAl2019}, and solving to global optimality using \cite{GurobiOptimization2020}, with a tolerance of $1\times 10^{-5}$. We solve all SDPs using Splitting Conic Solver (SCS) \cite{ODonoghueChuEtAl2016}, also with a tolerance of $1\times 10^{-5}$.

\begin{figure}
     \centering
    \begin{subfigure}[b]{0.49\columnwidth}
         \centering
         \includegraphics[width=\textwidth]{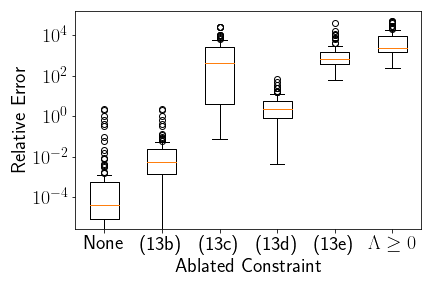}
         \caption{This figure plots the relative error when constraints in \eqref{opt:0SOS} are ablated. While removing \eqref{eq:input_lin} typically results in reasonable relaxation gaps, the removal of any other constraint renders the relaxation entirely uninformative with gaps typically on the order of $10^2$.}
    \label{fig:abl1}
     \end{subfigure}
     \hfill
     \begin{subfigure}[b]{0.49\columnwidth}
         \centering
         \includegraphics[width=\textwidth]{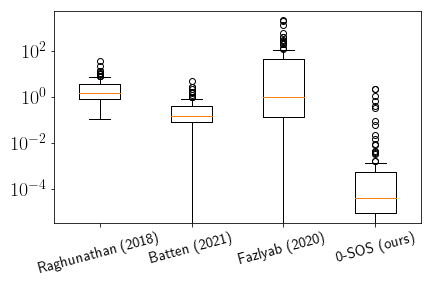}
             \caption{This figure shows the relative relaxation gap of the 0-SOS relaxation, \eqref{opt:0SOS}, as well as those of several existing relaxations. The 0-SOS relaxation is often exact, with a median error that is orders of magnitude smaller than the other relaxations.}
    \label{fig:comp1}
     \end{subfigure}
\vspace{-1.0cm}
\end{figure}

\paragraph{Ablating Verification-Defining Constraints.} In Section \ref{sec:results}, we have shown that removal of any of the verification-defining constraints fundamentally changes the underlying problem being relaxed, however we have left unanswered the degree to which removing individual constraints increases the relaxation gap.
We evaluate this empirically by considering variants of \eqref{opt:0SOS} where constraints are removed individually from the problem (termed the ablated SDPs). 
We focus on ablating \eqref{eq:input_lin}-\eqref{eq:nw_self_quad}, and $\Lambda \geq 0$, because constraints \eqref{eq:relu} and $\lambda \geq 0$, are not typically omitted in the literature. 
 The relative errors of each of the ablated SDPs are plotted in Figure \ref{fig:abl1}.

We find that removing the input linear constraints, \eqref{eq:input_lin}, has the smallest effect on the relaxation gap, even sometimes resulting in exact relaxations. On the other hand, removing the non-negativity constraint, $\Lambda \geq 0$, typically leads to bounds that are entirely uninformative, with relative errors ranging from 243 to $5.1\times10^4$. We also find that without the network self-quadratic constraints, \eqref{eq:nw_self_quad}, the relative errors range from 63 to $4\times 10^4$. The results of this ablation study are significant, because in the context of Figure \eqref{fig:chart}, we see that there are no existing approaches that entirely encode \eqref{eq:input_lin}, \eqref{eq:nw_self_quad}, and $\Lambda \geq 0$.

\paragraph{Comparison to the Literature.}
While we have shown that removing the verification-defining constraints typically introduces a significant relaxation gap, all existing work surveyed above also includes various strengthening constraints that imbue additional knowledge of the problem structure. 
In this experiment, we evaluate whether the strengthening constraints are sufficient to make up for the missing verification-defining constraints. 
Figure \ref{fig:comp1} plots the relative errors of the 0-SOS relaxation \eqref{opt:0SOS}, and that of several relaxations proposed in the literature.

We find that \eqref{opt:0SOS} typically results in the tightest relaxations, and is often exact. \cite{BattenKouvarosEtAl2021} and \cite{FazlyabMorariEtAl2022} are sometimes exact, while \cite{RaghunathanSteinhardtEtAl2018} never is. All exhibit median relaxation gaps that are orders of magnitudes larger than that of \eqref{opt:0SOS}. 
Notably, all exhibit considerably smaller relaxation gaps than the ablated SDPs.
This indicates that the gaps observed in the baselines are largely dependent on the strengthening constraints.

\section{CONCLUSION}\label{sec:conclusion}
In this paper, we have shown that verification of ReLU networks can be formulated exactly as a complete positive program, and provided analysis showing that the formulation is minimal. Leveraging the similarity between the form of the CPP formulation and SDP-based relaxations, we use our formulation as the basis for a unifying perspective on existing approaches. We also provide empirical evaluation demonstrating (1) that the verification defining constraints are indispensable for accurately representing the verification problem, and (2) the 0-SOS relaxation of the proposed formulation exhibits relaxation gaps that are orders of magnitude smaller than existing methods, and is often exact.

While we have focused on the 0-SOS relaxation as a comparison to existing approaches, it is just one in a hierarchy of relaxations now at our disposal.
Because higher order order SOS relaxations scale poorly, however, their application to solving verification requires further investigation.
Fortunately, there is an abundance of literature devoted to systematically improving the tractability of SOS-based methods. 
In particular, we highlight the r-DSOS and r-SDSOS hierarchies proposed by \cite{AhmadiMajumdar2019} as promising algorithmic paradigms for addressing the CPP formulation we have proposed.
We believe the most valuable contribution of this work is the potential for a new \emph{class} of verification methods that can systematically trade-off tightness and efficiency, with exactness as the limiting case.

\newpage
\textbf{Acknowledgements}

This work was supported by NSF CCF (grant \#1918549) and the NASA University Leadership Initiative (grant \#80NSSC20M0163); this article solely reflects the opinions and conclusions of its authors and not any NSF or NASA entity.

\bibliography{rabrown1}
\onecolumn
\appendix
\section{APPENDIX}
\renewcommand{\quadmat}[3]{\begin{pmatrix} #1 & #2 \\ #2 ^\top & #3 \end{pmatrix}}

\subsection{Illustration of expanded out matrices}
In this paper, we have made a number of routine conversions between matrices that act on pre/post-activation variables corresponding to a single layer, to those that act on the collection of all positive/negative splittings of variables from all layers. In this section, we provide a concrete example of how these conversions are carried out. Specifically, we illustrate the relationship between the matrices that act on single layers, $\overline{M}$, those that act on all variables, including slack variables $M$, and those that act on all variables except slack variables $\underline{M}$.

Consider an input set defined as $\mathcal{X} = \{x \in \mathbb{R}^{h_0} \mid \overline{A} x \leq a \}$. 
The matrix $\overline{A} \in \mathbb{R}^{m \times h_0}$ only acts on the input vector, $z_{0, *}$.
In this paper, we use the positive/negative splitting, so the input constraint is given by
\begin{equation}
    \overline{A}(\lambda^+_{0, *} - \lambda^-_{0, *})  \leq a.
\end{equation}
In addition, we concatenate all of the neurons and slack variables in a single vector
\begin{equation}
    \lambda 
    = 
    \begin{pmatrix}
        \lambda^+_{0, *}\\
        \lambda^+_{1, *}\\
        \vdots\\
        \lambda^+_{n, *}\\
        \lambda^-_{0, *}\\
        \lambda^-_{1, *}\\
        \vdots\\
        \lambda^-_{n, *}\\
        s
    \end{pmatrix}
\end{equation}
where $s \in \mathbb{R}^m$ is a vector of slack variables corresponding to the input inequality constraints. Let $N := \sum_{i = 0}^n h_i$ be the number of neurons, $\lambda \in \mathbb{R}^{2N + m}$ (the factor of two is because there are two variables for every neuron).
Then, the constraint $\overline{A}(\lambda^+_{0, *} - \lambda^-_{0, *})  \leq a$ may be written as
\begin{equation}
    \begin{pmatrix}
        \overline{A} & 0 & \hdots& 0 & -\overline{A} & 0 & \hdots & 0 & 0
    \end{pmatrix}
        \begin{pmatrix}
        \lambda^+_{0, *}\\
        \lambda^+_{1, *}\\
        \vdots\\
        \lambda^+_{n, *}\\
        \lambda^-_{0, *}\\
        \lambda^-_{1, *}\\
        \vdots\\
        \lambda^-_{n, *}\\
        s
    \end{pmatrix}
    \leq a.
\end{equation}
We define $\underline{A} \in \mathbb{R}^{m \times (2N + m)}$, which acts on all variables but does not incorporate the effect of slack variables, as 
\begin{equation}
\underline{A} := 
    \begin{pmatrix}
        \overline{A} & 0 & \hdots 0 & -\overline{A} & 0 & \hdots 0 & 0
    \end{pmatrix}.
\end{equation}
Similarly, for appropriate $s \geq 0$, then $\overline{A}(\lambda^+_{0, *} - \lambda^-_{0, *})  \leq a$ may be written as
\begin{equation}
    \begin{pmatrix}
        \overline{A} & 0 & \hdots& 0 & -\overline{A} & 0 & \hdots & 0 & I
    \end{pmatrix}
        \begin{pmatrix}
        \lambda^+_{0, *}\\
        \lambda^+_{1, *}\\
        \vdots\\
        \lambda^+_{n, *}\\
        \lambda^-_{0, *}\\
        \lambda^-_{1, *}\\
        \vdots\\
        \lambda^-_{n, *}\\
        s
    \end{pmatrix}
    = a.
\end{equation}
We define $A \in \mathbb{R}^{m \times (2N + m)}$, which acts on all variables including slacks, as 
\begin{equation}
    A :=
    \begin{pmatrix}
        \overline{A} & 0 & \hdots 0 & -\overline{A} & 0 & \hdots 0 & I
    \end{pmatrix}.
\end{equation}
\subsection{Proofs of Section \ref{sec:results}}

\begin{theorem*}[\ref{thm:verif_cpp}]
    Problem \eqref{opt:linear_exact} is equivalent to the following convex optimization problem:
    \begin{mini*}
    {\lambda,\, \Lambda}{c^\top \lambda}
    {}{\textup{\texttt{OPT}}_{\textup{\texttt{CPP}} = }}
    \addConstraint{A_{j, *} \lambda} {= a_j  \quad \tag{\ref{eq:input_lin}}}
    \addConstraint{\innerMat{A_{j, *}^\top A_{j, *}}{\Lambda}}{ = a^2_j\quad \tag{\ref{eq:input_self_quad}}}
    \addConstraint{W[i]_{j, *} \lambda }{= b[i]_j \quad \tag{\ref{eq:nw_lin}}}
    \addConstraint{\innerMat{W[i]_{j, *}^\top W[i]_{j, *}}{\Lambda}} {= b[i]^2_j \quad \tag{\ref{eq:nw_self_quad}}}
    \addConstraint{\Lambda[\lambda^+_{i, j}, \lambda^-_{i, j}]} {= 0 \quad \tag{\ref{eq:relu}}}
    \addConstraint{\quadmat{\Lambda}{\lambda}{1}}{\in \mathcal{C}^* \quad \tag{\ref{eq:comp_pos}}}
    \end{mini*}
    Exactness is defined by the following conditions: 
    \begin{enumerate}
        \item $\textup{\texttt{OPT}} =  \textup{\texttt{OPT}}_{\textup{\texttt{CPP}}}$: Problem \eqref{opt:linear_exact} and Problem \eqref{opt:CPP} have the same objective value.
        \item If $(\lambda^*, \Lambda^*) = \arg \min \eqref{opt:CPP}$ then $\lambda^*$ is in the convex hull of optimal solutions for Problem \eqref{opt:linear_exact}.
    \end{enumerate}
    \begin{proof}
    This follows as a direct application of Theorem 3.2 from \cite{Burer2009}. The two assumptions required are:
    \begin{itemize}
        \item The binary variables satisfy the \emph{key assumption}, chiefly, boundedness. This assumption is automatically satisfied as formulation \eqref{opt:linear_exact} does not have any binary variables.
        \item Each variable in a complementarity constraint is bounded. This follows from our assumption of a bounded input set. Explicit bounds on each neuron can be derived by forward propagating the input bounds.
    \end{itemize}
    
    \end{proof}
\end{theorem*}

\begin{theorem*}[\ref{thm:lin_eq_factor}]
    Suppose there is a factorization
    \begin{equation}
        \quadmat{\Lambda}{\lambda}{1} 
        = \sum_{k = 1}^K 
        \begin{pmatrix} \lambda^{(k)}\\ \xi^{(k)} \end{pmatrix}
        \begin{pmatrix} \lambda^{(k)}\\ \xi^{(k)} \end{pmatrix}^\top.
    \end{equation}
    Then, $ V_{i, *} \lambda^{(k)} = \xi^{(k)} v_i$ for all $k$ if and only if
    \begin{align}
        V_{i, *}\lambda &= v_i \tag{\ref{eq:lin_generic}}\\
        \innerMat{V_{i, *}^\top V_{i, *}}{\Lambda} &= v_i^2. \tag{\ref{eq:self_quad_generic}}
    \end{align}
    \begin{proof}
    $(\Leftarrow)$ We will first show that if $(\lambda, \, \Lambda)$ satisfy \eqref{eq:lin_generic} and \eqref{eq:self_quad_generic}, then $V_{i, *} \lambda^{(k)} = \xi^{(k)} v_i$.\\
    Expanding out the outer product, 
    \begin{equation}
        \quadmat{\Lambda}{\lambda}{1} 
        = \sum_{k = 1}^K 
        \begin{pmatrix} \lambda^{(k)}\\ \xi^{(k)} \end{pmatrix}
        \begin{pmatrix} \lambda^{(k)}\\ \xi^{(k)} \end{pmatrix}^\top
         = \sum_{k = 1}^K 
        \begin{pmatrix} \lambda^{(k)}(\lambda^{(k)})^\top & \xi^{(k)} \lambda^{(k)} \\
        \xi^{(k)} (\lambda^{(k)})^\top & (\xi^{(k)})^2 \end{pmatrix}.
    \end{equation} 
    Matching up terms, $\{\lambda^{(k)}, \, \xi^{(k)}\}$ must satisfy 
    \begin{align}
        \sum_k (\xi^{(k)})^2 &= 1\\
        \sum_k \xi^{(k)} \lambda^{(k)} &= \lambda\\
        \sum_k \lambda^{(k)} (\lambda^{(k)})^\top &= \Lambda
    \end{align}
    Substituting these summations into the constraints, we get that 
    \begin{align}
        &V_{i, *} \lambda = \sum_k \xi^{(k)} V_{i, *} \lambda^{(k)} = v_i \label{eq:lin_factor}\\
        &  \innerMat{V^\top_{i, *} V_{i, *}}{\Lambda} = \sum_k (V_{i, *} \lambda^{(k)})^2 = v_i^2
    \end{align}
    By squaring \eqref{eq:lin_factor} and using the fact that $\sum_k (\xi^{(k)})^2 = 1$, we get 
    \begin{align}
        v_i^2 &= (V_{i, *} \lambda)^2 \\
        &= (\sum_k \xi^{(k)} V_{i, *} \lambda^{(k)})^2 \\
        &\leq (\sum_k (\xi^{(k)})^2) \sum_k (V_{i, *} \lambda^{(k)})^2 \label{eq:cs}\\
        &= \sum_k (V_{i, *}\lambda^{(k)})^2 \\
        &= v_i^2
    \end{align}
    where \eqref{eq:cs} follows from Cauchy-Schwarz. More specifically, this shows that in \eqref{eq:cs} Cauchy-Schwarz holds with equality, which means that $[\xi^{(k)}]_k$ and $[V_{i, *} \lambda^{(k)}]_k$ are collinear. This means that there is a scalar $\alpha \in \mathbb{R}$ such that $\alpha \xi^{(k)} = V_{i, *}\lambda^{(k)}$ for all $k$. Substituting this into \eqref{eq:lin_factor}, we get 
    \begin{equation}
        v_i = \sum_k \xi^{(k)} V_{i, *} \lambda^{(k)} = \sum_k (\xi^{(k)})^2 \alpha = \alpha.
    \end{equation}
    This shows $V_{i, *}\lambda^{(k)} = \xi^{(k)} v_i$, thus concluding the proof of the reverse implication.\\
    $(\Rightarrow)$ We now just have to show that if $V_{i, *}\lambda^{(k)} = \xi^{(k)} v_i$ for all $k$ then $V_{i, *}\lambda = v_i$ and $\innerMat{V^\top_{i, *} V_{i, *}}{\Lambda} = v_i^2$.
    The former follows by substituting in the relation $\sum_k \xi^{(k)} \lambda^{(k)} = \lambda$
    \begin{align}
        V_{i, *} \lambda &= \sum_k \xi^{(k)} V_{i, *} \lambda^{(k)}\\
        &= \sum_k (\xi^{(k)})^2 v_i\\
        &= v_i
    \end{align}
    where the last line follows because $\sum_k (\xi^{(k)})^2 = 1$.
    The latter follows from substituting in the relationship $\innerMat{V^\top_{i, *} V_{i, *}}{\Lambda} = \sum_k (V_{i, *}\lambda^{(k)})^2 $:
    \begin{align}
        \innerMat{V^\top_{i, *} V_{i, *}}{\Lambda} &= \sum_k \innerMat{V^\top_{i, *} V_{i, *}}{\lambda^k (\lambda^k)^\top}\\
        &= \sum_k (V_{i, *} \lambda^k)^2 \\
        &= \sum_k (\xi^k)^2 v_i^2\\
        &= v_i^2.
    \end{align}
    \end{proof}
\end{theorem*}

\begin{corollary*} [\ref{cor:ineq_cpp_0sos}]
   Suppose there is a factorization,
    \begin{equation} \tag{\ref{eq:slack_factor}}{
        \begin{pmatrix}
        \Lambda[\lambda, \lambda^\top] & \Lambda[\lambda, s^\top] & \lambda\\
        \Lambda[s, \lambda ^\top] & \Lambda[s, s^\top] & s\\
        \lambda^\top & s^\top & 1
        \end{pmatrix}
         = \sum_{k = 1}^K 
        \begin{pmatrix} \lambda^{(k)}\\ s^{(k)} \\\xi^{(k)} \end{pmatrix}
        \begin{pmatrix} \lambda^{(k)}\\ s^{(k)} \\\xi^{(k)} \end{pmatrix}^\top}
    \end{equation}
    
    and for all $i = \{1, \ldots L\}$, $(\lambda, \, \Lambda)$ satisfy
    \begin{align}
        &V_{i, *}\lambda + s_i = v_i, \tag{\ref{eq:slack_lin}}\\
        &\innerMat{V^\top_{i, *} V_{i, *}}{\Lambda[\lambda, \lambda^\top]} + 2 V_{i, *} \Lambda[\lambda, s_i] + \Lambda[s_i, s_i] = v_i^2. \tag{\ref{eq:slack_quad}}
    \end{align}

    Then we must have:
    \begin{enumerate}
        \item For $\begin{pmatrix}
        \Lambda[\lambda, \lambda^\top] & \Lambda[\lambda, s^\top] & \lambda\\
        \Lambda[s, \lambda ^\top] & \Lambda[s, s^\top] & s\\
        \lambda^\top & s^\top & 1
        \end{pmatrix} \in C^*$,
        each factor individually satisfies the inequalities:
        \begin{equation}{
            V_{i, *}\lambda^{(k)} \leq \xi^{(k)} v_i.}
        \end{equation}
        \item For $\begin{pmatrix}
        \Lambda[\lambda, \lambda^\top] & \Lambda[\lambda, s^\top] & \lambda\\
        \Lambda[s, \lambda ^\top] & \Lambda[s, s^\top] & s\\
        \lambda^\top & s^\top & 1
        \end{pmatrix} \in S^+$, and $s \geq 0$
        the weighted sum of factors satisfies the inequalities:
        \begin{equation}{
            V_{i, *} \lambda = \sum_k \xi^{(k)} V_{i, *}  \lambda^{(k)} \leq v_i}.
        \end{equation}
    \end{enumerate}
    
    \begin{proof}
    Leveraging Theorem \ref{thm:lin_eq_factor}, we know that if 
    \begin{equation}\tag{\ref{eq:slack_factor}}
        \begin{pmatrix}
        \Lambda[\lambda, \lambda^\top] & \Lambda[\lambda, s^\top] & \lambda\\
        \Lambda[s, \lambda ^\top] & \Lambda[s, s^\top] & s\\
        \lambda^\top & s^\top & 1
        \end{pmatrix}
         = \sum_{k = 1}^K 
        \begin{pmatrix} \lambda^{(k)}\\ s^{(k)} \\\xi^{(k)} \end{pmatrix}
        \begin{pmatrix} \lambda^{(k)}\\ s^{(k)} \\\xi^{(k)} \end{pmatrix}^\top
    \end{equation}
    and $(\lambda, \Lambda)$ satisfy \eqref{eq:slack_lin} and \eqref{eq:slack_quad}, then $V_{i, *}\lambda^{(k)} + s^{(k)}_i = \xi^{(k)} v_i$ for all $k$. 
    
    If 
    $\begin{pmatrix}
        \Lambda[\lambda, \lambda^\top] & \Lambda[\lambda, s^\top] & \lambda\\
        \Lambda[s, \lambda ^\top] & \Lambda[s, s^\top] & s\\
        \lambda^\top & s^\top & 1
        \end{pmatrix} \in C^*$ then in particular $s^{(k)} \geq 0$ so $V_{i, *} \lambda^{(k)} \leq \xi^{(k)} v_i$.
        
    In the case where 
        $\begin{pmatrix}
        \Lambda[\lambda, \lambda^\top] & \Lambda[\lambda, s^\top] & \lambda\\
        \Lambda[s, \lambda ^\top] & \Lambda[s, s^\top] & s\\
        \lambda^\top & s^\top & 1
        \end{pmatrix} \in S^+$, and $s \geq 0$, the result follows directly from the constraints $V_{i, *} \lambda + s_i = v_i $ and $s_i \geq 0$.
    \end{proof}
\end{corollary*}

\begin{corollary*}[\ref{cor:cross_redundant}]
    If $\begin{pmatrix}
        \Lambda[\lambda, \lambda^\top] & \Lambda[\lambda, s^\top] & \lambda\\
        \Lambda[s, \lambda ^\top] & \Lambda[s, s^\top] & s\\
        \lambda^\top & s^\top & 1
        \end{pmatrix} \in S^+ \cap \mathcal{N}$, and the following hold:
        \begin{enumerate}
            \item There is a factorization 
            \begin{equation}\tag{\ref{eq:slack_factor}}
                \begin{pmatrix}
                    \Lambda[\lambda, \lambda^\top] & \Lambda[\lambda, s^\top] & \lambda\\
                    \Lambda[s, \lambda ^\top] & \Lambda[s, s^\top] & s\\
                    \lambda^\top & s^\top & 1
                \end{pmatrix}
                 = \sum_{k = 1}^K 
                \begin{pmatrix} \lambda^{(k)}\\ s^{(k)} \\\xi^{(k)} \end{pmatrix}
                \begin{pmatrix} \lambda^{(k)}\\ s^{(k)} \\\xi^{(k)} \end{pmatrix}^\top
            \end{equation}
            
            \item     
            \begin{equation}\tag{\ref{eq:slack_lin}}
                V_{i, *}\lambda + s_i = v_i \end{equation}
            \item 
            \begin{equation}\tag{\ref{eq:slack_quad}}
              \innerMat{V^\top_{i, *} V_{i, *}}{\Lambda[\lambda, \lambda^\top]} + 2 V_{i, *} \Lambda[\lambda, s_i] + \Lambda[s_i, s_i] = v_i^2  
            \end{equation}
        \end{enumerate}
            then the following inequality also holds:   
        \begin{equation}\tag{\ref{eq:cross_quad_constr}}{
            v_i v_j - v_i V_{j, *} \lambda  - v_j V_{i, *}\lambda  + \innerMat{V^\top_{i, *} V_{j, *}}{\Lambda[\lambda, \lambda^\top]} \geq 0}
        \end{equation}
    \begin{proof}
    
    From Theorem \ref{thm:lin_eq_factor}, $V_{i, *} \lambda^{(k)} + s^{(k)}_i = \xi^{(k)} v_i$ holds.
    Substituting we have 
    \begin{align}
        & v_iv_j - v_i V_{j, *} \lambda  - v_j V_{i, *}\lambda  + \innerMat{V_{i,*}^\top V_{j, *}}{\Lambda[\lambda, \lambda^\top]} \\
        &=v_i v_j + \sum_k  \left(-v_i V_{j, *} (\xi^{(k)}\lambda^{(k)})  -v_j V_{i, *} (\xi^{(k)}\lambda^{(k)}) + \innerMat{V_{i,*}^\top V_{j, *}}{\lambda^{(k)} (\lambda^{(k)})^\top}\right)\\
        &= \sum_k  (\xi^{(k)})^2 v_i v_j - \xi^{(k)} v_i V_{j, *} \lambda^{(k)}  - \xi^{(k)}v_j V_{i, *} \lambda^{(k)} + V_{i,*}\lambda^{(k)} V_{j, *}\lambda^{(k)}\\
        &= \sum_k (\xi^{(k)}v_i - V_{i, *}\lambda^{(k)})(\xi^{(k)}v_j - V_{j, *}\lambda^{(k)})\\
        &= \sum_k s^{(k)}_i s^{(k)}_j\\
        &\geq 0.
    \end{align}
    The last line follows from the equality 
    \begin{equation}
        \Lambda[s_i, s_j] = \sum_k s^{(k)}_i s^{(k)}_j
    \end{equation}
    and the entrywise non-negativity of $\begin{pmatrix}
        \Lambda[\lambda, \lambda^\top] & \Lambda[\lambda, s^\top] & \lambda\\
        \Lambda[s, \lambda ^\top] & \Lambda[s, s^\top] & s\\
        \lambda^\top & s^\top & 1
        \end{pmatrix} \in \mathcal{N}$.
    \end{proof}
    \end{corollary*}

    \subsection{Reformulation of \cite{RaghunathanSteinhardtEtAl2018}}
    \cite{RaghunathanSteinhardtEtAl2018} proposed the following SDP relaxation:
    \begin{mini*}|s|
    {z, Z}{\overline{c}^\top z_{n, *}}
    {}{\tag{\ref{opt:raghunation_raw}}}
    \addConstraint{z_{i, *} \geq 0, \, z_{i, *} \geq \overline{W}[i - 1] z_{i - 1, *}}
    \addConstraint{Z[z_{i, j}, z_{i, j}] = \overline{W}[i - 1]_{j, *} Z[z_{i - 1, *}, z_{i,j}]}
    \addConstraint{-Z[z_{i,j}, z_{i,j}] + (u_{i, j} + l_{i, j})z_{i,j} - l_{i, j}u_{i,j} \geq 0, \, \forall i, j}
    \addConstraint{\quadmat{Z}{z}{1} \in \mathcal{S}^+}
    \end{mini*}
    
    The proposed SDP relaxation is based on directly transcribing the constraints in the following quadratically constrained quadratic program:
    \begin{mini}|s|
    {z}{\overline{c}^\top z_{n, *}}
    {}{\label{opt:raghunathan_qcqp}}
    \addConstraint{z_{i, *} \geq 0, \, z_{i, *} \geq \overline{W}[i - 1] z_{i - 1, *}}
    \addConstraint{z_{i,j}^2 = z_{i,j} \overline{W}[i - 1]_{j, *}z_{i - 1, *} }
    \addConstraint{(z_{i,j} - l_{i, j})(u_{i, j} - z_{i, j}) \geq 0, \, \forall i, j}
    \end{mini}

    To clarify the relationship with our proposed framework, we include the biases, introduce an additional equality constraint $\hat{z}_{i+1, *} = W[i] z_{i, *} + b[i]$ and the change the inequality from $z_{i+1, *} \geq W[i] z_{i, *} + b[i]$ to $z_{i, *} \geq \hat{z}_{i, *}$. Without loss of generality, we will assume that $\hat{l}_{i, j} < 0 < \hat{u}_{i, j}$ for all $i \geq 1$, otherwise the neuron $(i, j)$ could be treated as the identity, if $\hat{l}_{i, j} \geq 0$, or zero, if $\hat{u}_{i, j}\leq 0$. The post-activation bounds are then $l_{i, j} = 0$ and $u_{i, j} = \hat{u}_{i, j}$, and the constraint $(z_{i,j} - l_{i, j})(u_{i, j} - z_{i, j}) \geq 0$ becomes $z_{i,j} (\hat{u}_{i, j} - z_{i, j}) \geq 0$. With a change of variables to the positive/negative splitting, \eqref{opt:raghunathan_qcqp} is equivalent to:
    \begin{mini}|s|
    {\lambda}{\overline{c}^\top (\lambda^+_{n, *} - \lambda^-_{n, *})}
    {}{\label{opt:raghunathan_qcqp_pn}}
    \addConstraint{\lambda^+_{i, *} \geq 0, \, \lambda^-_{i, *} \geq 0, \, \lambda^+_{i + 1, *} - \lambda^-_{i+1, *} = \overline{W}[i] \lambda^+_{i, *} + b[i]}
    \addConstraint{\lambda^+_{i, *} \lambda^-_{i, *} = 0}
    \addConstraint{\lambda^+_{i,j} (\hat{u}_{i, j} - \lambda^+_{i, j}) \geq 0, \, i \geq 1}
    \addConstraint{((\lambda^+_{0, *} - \lambda^-_{0, *}) - l_{0, *})(u_{0, *} - (\lambda^+_{0, *} - \lambda^-_{0, *}) ) \geq 0}
    \end{mini}

    Directly transcribing the constraints of \eqref{opt:raghunathan_qcqp_pn} results in the following SDP:
    \begin{mini*}|s|
    {\lambda, \, \Lambda}{c^\top \lambda}
    {}{\tag{\ref{opt:SDR}}}
    \addConstraint{(l_{0, i} + u_{0, i}) U_{i, *}\lambda - \innerMat{U^\top_{i, *}U_{i, *}}{\Lambda} - u_il_i \geq 0}
    \addConstraint{\hat{u}_{i, j} \lambda^+_{i, j} -  \Lambda[\lambda^+_{i, j}, \lambda^+_{i, j}] \geq 0, \, i \geq 1}
    \addConstraint{W[i]_{j, *} \lambda = b[i]_j}
    \addConstraint{\Lambda[\lambda^+_{i, j}, \lambda^-_{i, j}] = 0}
    \addConstraint{\quadmat{\Lambda}{\lambda}{1} \in S^+}
    \addConstraint{\lambda \geq 0}
    \end{mini*}
    where $U$ is defined as $U_{i, *} \lambda = \lambda^+_{0, i} - \lambda^-_{0, i}$.
    
    \subsubsection{Extensions to \cite{RaghunathanSteinhardtEtAl2018}}
    \paragraph{Direct Extensions.}
    The following works have focused on improving some of the computational aspects of \cite{RaghunathanSteinhardtEtAl2018} (e.g., time to solution, memory usage) without aiming to improve the relaxation gap of \cite{RaghunathanSteinhardtEtAl2018}.
    
    The construction in \cite{DvijothamStanforthEtAl2020} is based on the diagonally dominant SOS hierarchy (0-DSOS) relaxation. It approximates $S^+$ with the cone of diagonally dominant matrices, ultimately resulting in an LP. This is a further relaxation of \eqref{opt:SDR}, and reduces computation time at the expense of a larger relaxation gap.
    
    Based on the fact that memory usage, rather than compute, is the main barrier for scaling of SDPs, \cite{DathathriDvijothamEtAl2020} proposes solving \eqref{opt:SDR} using first order methods. 
    This comes at the expense of increased time to solution.

    \paragraph{SDP-based Branching.}
    Previously, \cite{Zhang2020} defined an SDP relaxation to be tight if it has a unique rank-one solution. Based on this, \cite{AndersonMaEtAl2021} proposed a metric measuring how far the solution is from being rank one and shows that a uniform partitioning of the input (before neuron bound propagation) leads to the greatest reduction of an upper bound of this metric. This approach can be thought of as branching on the input.

    Critically, Theorem \ref{thm:verif_cpp} states that the CPP \eqref{opt:CPP} is always exact, even if the optimal solution is not rank one. 
    If an optimizer of the (0-SOS) relaxation \eqref{opt:0SOS} happens to be completely positive, it is exact even if it is not rank one.
    This means that the rank-one condition is a sufficient but not necessary condition for exactness.
    For this reason, we propose a less-restrictive definition of tightness, allowing exact but non-rank-one optimizers; these points are further clarified in Section \ref{subsec:rank}.

    \cite{MaSojoudi2020}  work directly with \eqref{opt:SDR}, and propose a method equivalent to spatial branch and bound. In particular, they choose a basis $\{\phi_i\}$ for $\mathbb{R}^N$, and partition the search space along the axes aligned with $\{\phi_i\}$. To illustrate how their approach works, consider $\{\phi_i\}$ given by the standard basis vectors, and a partitioning of $\lambda^+_{i, j}$. Their approach generates $M$ sub-problems that are determined by points $l_{i, j} = \gamma_0 < \gamma_1 < \ldots < \gamma_M = u_{i, j}$, where the $m$th sub-problem is defined by the following additional constraints:
        \begin{align}
            &\lambda^+_{i, j} \leq \gamma_m\\
            &\lambda^+_{i, j} \geq \gamma_{m - 1}\\
            &\gamma_m \lambda^+_{i, j} - \gamma_m \gamma_{m - 1} - \Lambda[\lambda^+_{i, j}, \lambda^+_{i, j}] + \gamma_m \lambda^+_{i, j} \geq 0 \label{eq:ma_cross_quad}
        \end{align}
        Notice that equation \eqref{eq:ma_cross_quad} is the cross-quadratic constraint derived from $\lambda^+_{i, j}  \leq \gamma_m$ and $\lambda^+_{i, j}\geq \gamma_{m - 1}$.

    \subsection{Reformulation of \cite{FazlyabMorariEtAl2022}}
    \cite{FazlyabMorariEtAl2022} develop a method where they abstract the definition of particular constraint sets as quadratic constraints. In contrast to the quadratic constraints presented thus far in this paper, the quadratic constraints in \cite{FazlyabMorariEtAl2022} are defined by (potentially infinite) matrix cones.
    
    \subsubsection{Background on quadratic constraints}
    For example, if $\mathcal{X}$ is the input set, then the matrix cone describing $\mathcal{X}$ is defined as
    
    \begin{equation}
        \mathcal{P}_{\mathcal{X}} := \left\{P \mid 
        \begin{pmatrix} \lambda\\ 1\end{pmatrix}^\top 
        P
        \begin{pmatrix} \lambda \\ 1\end{pmatrix}
        \geq 0, \, \forall \lambda \in \mathcal{X}
        \right\}
    \end{equation}
    
    In the converse, $\mathcal{P}_{\mathcal{X}}$ can be used to over-approximate $\mathcal{X}$ via an infinite set of constraints:
    \begin{equation}
        \mathcal{X} \subseteq \bigcap_{P \in \mathcal{P}_{\mathcal{X}}}
        \left\{
        \lambda \mid 
        \begin{pmatrix} \lambda \\ 1\end{pmatrix}^\top 
        P
        \begin{pmatrix} \lambda \\ 1\end{pmatrix}
        \geq 0
        \right\}
    \end{equation}
    
    The quadratic constraints are related to the SDP/CPP formulations of verification through the following equality:
    \begin{equation}
        \begin{pmatrix} \lambda \\ 1\end{pmatrix}^\top 
        P
        \begin{pmatrix} \lambda \\ 1\end{pmatrix}
        = 
        \innerMat{P}{\begin{pmatrix} \lambda \\ 1\end{pmatrix} \begin{pmatrix} \lambda \\ 1\end{pmatrix}^\top}.
    \end{equation}
    The premise of the SDP/CPP formulations is relaxing the rank-one outer product $\begin{pmatrix} \lambda \\ 1\end{pmatrix} \begin{pmatrix} \lambda \\ 1\end{pmatrix}^\top$ to $\quadmat{\Lambda}{\lambda}{1} \in \mathcal{S}^+$.
    
    \subsubsection{Equivalence between atomic and infinite quadratic constraints}
    Our analysis is based on identifying an equivalent, finite set of quadratic constraints (``atomic constraints") for each infinite set of quadratic constraints proposed in \cite{FazlyabMorariEtAl2022}---these correspond with the extreme rays of $\mathcal{P}$. In other words, for each quadratic constraint set $\mathcal{P}$ with infinite cardinality, we will find a finite set $\hat{\mathcal{P}}$ such that
    \begin{equation}
        \bigcap_{P \in \mathcal{P}}
        \left\{
        \lambda \mid 
        \begin{pmatrix} \lambda \\ 1\end{pmatrix}^\top 
        P
        \begin{pmatrix} \lambda \\ 1\end{pmatrix}
        \geq 0
        \right\}
        =
        \bigcap_{P \in \hat{\mathcal{P}}}
        \left\{
        \lambda \mid 
        \begin{pmatrix} \lambda \\ 1\end{pmatrix}^\top 
        P
        \begin{pmatrix} \lambda \\ 1\end{pmatrix}
        \geq 0
        \right\}.
    \end{equation}
    
    In particular, if 
    \begin{equation}
      \mathcal{P} = \left\{\sum_{i = 1}^{N_{\text{ineq}}} \alpha_i P_{\text{ineq}}^{(i)} + \sum_{i = 1}^{N_{\text{eq}}} \beta_i P_{\text{eq}}^{(i)} \mid \alpha_i \in \mathbb{R}_{\scriptscriptstyle \geq 0}, \, \beta_i \in \mathbb{R}
      \right \}
    \end{equation}
    then
    \begin{equation}
    \begin{aligned}
    \bigcap_{P \in \mathcal{P}}
        \left\{
        \lambda \mid 
        \begin{pmatrix} \lambda \\ 1\end{pmatrix}^\top 
        P
        \begin{pmatrix} \lambda \\ 1\end{pmatrix}
        \geq 0
        \right\}
    = 
    \Bigg \{
    \lambda \mid
    &\begin{pmatrix} \lambda \\ 1\end{pmatrix}^\top 
        P_{\text{ineq}}^{(i)}
    \begin{pmatrix} \lambda \\ 1\end{pmatrix}
    \geq 0
    , i = 1, \ldots, N_{\text{ineq}}, \\
    &\begin{pmatrix} \lambda \\ 1\end{pmatrix}^\top 
        P_{\text{eq}}^{(j)}
    \begin{pmatrix} \lambda \\ 1\end{pmatrix}
    = 0
    , j = 1, \ldots, N_{\text{eq}}
    \Bigg\}.
    \end{aligned}
    \end{equation}
    The equivalence can be shown as follows:
    \begin{itemize}
        \item If $\lambda$ satisfies
        \begin{equation}
            \begin{aligned}
                &\begin{pmatrix} \lambda \\ 1\end{pmatrix}^\top 
                    P_{\text{ineq}}^{(i)}
                \begin{pmatrix} \lambda \\ 1\end{pmatrix}
                \geq 0
                , i = 1, \ldots, N_{\text{ineq}}, \\
                &\begin{pmatrix} \lambda \\ 1\end{pmatrix}^\top 
                    P_{\text{eq}}^{(j)}
                \begin{pmatrix} \lambda \\ 1\end{pmatrix}
                = 0
                , j = 1, \ldots, N_{\text{eq}},
            \end{aligned}
        \end{equation}
        then for any $P \in \mathcal{P}$ using the decomposition $P = \sum_{i = 1}^{N_{\text{ineq}}} \alpha_i P_{\text{ineq}}^{(i)} + \sum_{i = 1}^{N_{\text{eq}}} \beta_i P_{\text{eq}}^{(i)}$ with $\alpha_i \geq 0$, 
            \begin{align}
                \begin{pmatrix} \lambda \\ 1\end{pmatrix}^\top 
                    P
                \begin{pmatrix} \lambda \\ 1\end{pmatrix}
                &=
                \begin{pmatrix} \lambda \\ 1\end{pmatrix}^\top 
                \left (\sum_{i = 1}^{N_{\text{ineq}}} \alpha_i P_{\text{ineq}}^{(i)} + \sum_{i = 1}^{N_{\text{eq}}} \beta_i P_{\text{eq}}^{(i)} \right)
                \begin{pmatrix} \lambda \\ 1\end{pmatrix} \\
                &= 
                \sum_{i = 1}^{N_{\text{ineq}}} \alpha_i \begin{pmatrix} \lambda \\ 1\end{pmatrix}^\top P_{\text{ineq}}^{(i)}\begin{pmatrix} \lambda \\ 1\end{pmatrix} + \sum_{i = 1}^{N_{\text{eq}}} \beta_i \begin{pmatrix} \lambda \\ 1\end{pmatrix}^\top P_{\text{eq}}^{(i)} 
                \begin{pmatrix} \lambda \\ 1\end{pmatrix} \\
                &\geq 0.
            \end{align}
            So 
            \begin{equation}
                \begin{aligned}
                \bigcap_{P \in \mathcal{P}}
                    \left\{
                    \lambda \mid 
                    \begin{pmatrix} \lambda \\ 1\end{pmatrix}^\top 
                    P
                    \begin{pmatrix} \lambda \\ 1\end{pmatrix}
                    \geq 0
                    \right\}
                \supseteq
                \Bigg \{
                \lambda \mid
                &\begin{pmatrix} \lambda \\ 1\end{pmatrix}^\top 
                    P_{\text{ineq}}^{(i)}
                \begin{pmatrix} \lambda \\ 1\end{pmatrix}
                \geq 0
                , i = 1, \ldots, N_{\text{ineq}}, \\
                &\begin{pmatrix} \lambda \\ 1\end{pmatrix}^\top 
                    P_{\text{eq}}^{(j)}
                \begin{pmatrix} \lambda \\ 1\end{pmatrix}
                = 0
                , j = 1, \ldots, N_{\text{eq}}
                \Bigg\}.
                \end{aligned}
            \end{equation}
        \item 
        Given $\lambda$, suppose there exists $i^*$ such that $\begin{pmatrix} \lambda \\ 1\end{pmatrix}^\top  P_{\text{ineq}}^{(i^*)} 
        \begin{pmatrix} \lambda \\ 1\end{pmatrix} 
        < 0$. Then, letting $\alpha_{i^*} = 1$, $\alpha_{i} = 0$ for $i \neq 0$ and $\beta_j = 0$, gives an example of $P=  \sum_{i = 1}^{N_{\text{ineq}}} \alpha_i P_{\text{ineq}}^{(i)} + \sum_{i = 1}^{N_{\text{eq}}} \beta_i P_{\text{eq}}^{(i)} \in \mathcal{P}$ with 
        \begin{equation}
        \begin{pmatrix} \lambda \\ 1\end{pmatrix}^\top 
            P
        \begin{pmatrix} \lambda \\ 1\end{pmatrix}
        < 0 \end{equation}
        An analogous argument holds if  there exists $i^*$ such that $\begin{pmatrix} \lambda \\ 1\end{pmatrix}^\top  P_{\text{eq}}^{(i^*)} 
        \begin{pmatrix} \lambda \\ 1\end{pmatrix} 
        \neq 0$ with $\beta_{i^*} = -\text{sign}\left(\begin{pmatrix} \lambda \\ 1\end{pmatrix}^\top  P_{\text{eq}}^{(i^*)} 
        \begin{pmatrix} \lambda \\ 1\end{pmatrix}\right)$.
        
        This shows that 
                    \begin{equation}
                \begin{aligned}
                \bigcap_{P \in \mathcal{P}}
                    \left\{
                    \lambda \mid 
                    \begin{pmatrix} \lambda \\ 1\end{pmatrix}^\top 
                    P
                    \begin{pmatrix} \lambda \\ 1\end{pmatrix}
                    \geq 0
                    \right\}
                \subseteq
                \Bigg \{
                \lambda \mid
                &\begin{pmatrix} \lambda \\ 1\end{pmatrix}^\top 
                    P_{\text{ineq}}^{(i)}
                \begin{pmatrix} \lambda \\ 1\end{pmatrix}
                \geq 0
                , i = 1, \ldots, N_{\text{ineq}}, \\
                &\begin{pmatrix} \lambda \\ 1\end{pmatrix}^\top 
                    P_{\text{eq}}^{(j)}
                \begin{pmatrix} \lambda \\ 1\end{pmatrix}
                = 0
                , j = 1, \ldots, N_{\text{eq}}
                \Bigg\}.
                \end{aligned}
            \end{equation}
    \end{itemize}
    Combining the inclusions in both directions shows that the two constraints are equal.
    \subsubsection{Polytopic input constraints}
    To encode polytopic input sets, $\mathcal{X} =\{x \in \mathbb{R}^{k_0} \mid \overline{A} x \leq a\}$, \cite{FazlyabMorariEtAl2022} use the following quadratic constraint set:
    \begin{equation}
        \mathcal{P}_{\mathcal{X}} = \left\{ P \mid P = 
        \begin{bmatrix}
        \underline{A}^\top \Gamma \underline{A} & -\underline{A}^\top \Gamma a\\
        -a^\top \Gamma \underline{A} & a^\top \Gamma a
        \end{bmatrix}
        \right\}
    \end{equation}
    where $\Gamma \geq 0$, $\Gamma_{i,i} = 0$ and $\Gamma = \Gamma^\top$. Each $P \in \mathcal{P}_{\mathcal{X}}$ can be decomposed as
    \begin{equation}
        P = \sum_{i \neq j} \Gamma_{i, j}
        \begin{bmatrix}
        \underline{A}^\top_{i, *} \underline{A}_{j, *} & \frac{-a_j \underline{A}^\top_{i, *} - a_i \underline{A}^\top_{j, *} }{2}\\
        \frac{-a_j \underline{A}_{i, *} - a_i \underline{A}_{j, *} }{2} & a_i a_j\\
        \end{bmatrix}
    \end{equation}
    As a consequence, $\mathcal{P}_{\mathcal{X}}$ defines constraints of the form:
    \begin{equation}
        \sum_{i \neq j} \Gamma_{i, j} (\underline{A}_{i, *} \lambda  - a_i)(\underline{A}_{j, *}\lambda - a_j) \geq 0
    \end{equation}
    with parameters $\Gamma_{i, j} \geq 0, \, \Gamma_{i, i} = 0$. 
    The equivalent atomic constraints are
    \begin{equation}
        (\underline{A}_{i, *} \lambda  - a_i)(\underline{A}_{j, *}\lambda - a_j) \geq 0, \quad \forall i \neq j.
    \end{equation}

    \subsubsection{ReLU constraints}
    \cite{FazlyabMorariEtAl2022} defined a global quadratic constraint for ReLUs, acting on $\begin{pmatrix} z \\ \hat{z} \\ 1\end{pmatrix}$, through a constraint set of the form 
    \begin{equation}
        \mathcal{P}
        = 
        \begin{bmatrix}
         P_{1, 1} & P_{1, 2} & P_{1, 3}\\
         P_{2, 1} & P_{2, 2} & P_{2, 3}\\
         P_{3, 1} & P_{3, 2} & P_{3, 3}
        \end{bmatrix}
    \end{equation}
    where 
    \begin{align}
    P_{1, 1} &= 0  \in \mathbb{R}^{N\times N}\\
    P_{1, 2} &= \text{diag}(\rho) + T \in \mathbb{R}^{N \times N}\\
    P_{1, 3} &= -\nu \in \mathbb{R}^{N}\\
    P_{2, 2} &= -2(\text{diag}(\rho) + T) \in \mathbb{R}^{N \times N}\\ 
    P_{2, 3} &= \nu + \eta \in \mathbb{R}^{N}\\
    P_{3, 3} &= 0 \in \mathbb{R}
    \end{align}
    where $\eta, \nu \geq 0$, $\rho$ is unconstrained, and the matrix $T$ is defined as 
    \begin{equation}
        T := \sum_{1 \leq i < j \leq N} \gamma_{i, j} (e_i - e_j)
        (e_i - e_j)^\top 
    \end{equation}
    where $e_i$ is the ith basis vector, and  $\gamma_{i, j} \geq 0$.
    Each quadratic constraint, $P$, acting on  $\begin{pmatrix} z \\ \hat{z} \\ 1\end{pmatrix}$ can be converted to an equivalent quadratic constraint acting on $\begin{pmatrix} \lambda^+ \\ \lambda^- \\ 1\end{pmatrix}$ by pre- and post-multiplying by $\begin{pmatrix} I & 0 & 0 \\ I & -I & 0 \\ 0 & 0 & 1 \end{pmatrix}$
    
    We have taken the liberty to substitute in $\alpha = 0, \, \beta = 1$ in $z_{i, j} = \max(\alpha \hat{z}_{i, j},\, \beta\hat{z}_{i, j})$. We have also renamed constraints to prevent notational clash with our notation.
    
    The quadratic constraint set, $\mathcal{P}$, thus enforces quadratic constraints of the form \cite[Lemma 3 and Appendix D]{FazlyabMorariEtAl2022}:
    \begin{equation}
        \begin{aligned}
        0 \geq  &\sum_{(i, j)} \left[\rho_{(i, j)} \lambda^+_{i, j}\lambda^-_{i, j} - \nu_{i, j}  \lambda^+_{i, j}  - \eta_{i, j} \lambda^-_{i, j} \right]\\
        &+ \sum_{(i, j) \neq (k, l)} \gamma_{(i, j), (k, l)}(\lambda^+_{i, j} - \lambda^+_{k, l})(\lambda^-_{i, j} - \lambda^-_{k, l})
        \end{aligned}
    \end{equation}
    Each term in the sum decomposes into the atomic constraints (together, equivalent to the full sum)
    \begin{align}
        \lambda^+_{i, j} &\geq 0 \\
        \lambda^-_{i, j} &\geq 0 \\ 
        \lambda^+_{i, j} \lambda^-_{i, j} &= 0 \\
        (\lambda^+_{i, j} - \lambda^+_{k, l})(\lambda^-_{i, j} - \lambda^-_{k, l}) &\leq 0 \label{eq:cross_relu}
    \end{align}
    Expanding out \eqref{eq:cross_relu}, 
    \begin{align}
        0 &\geq (\lambda^+_{i, j} - \lambda^+_{k, l})(\lambda^-_{i, j} - \lambda^-_{k, l})\\
        &= \lambda^+_{i, j}\lambda^-_{i, j} - \lambda^+_{i, j}\lambda^-_{k, l}  - \lambda^+_{k, l}\lambda^-_{i, j} + \lambda^+_{k, l}\lambda^-_{k, l}\\
        &= - \lambda^+_{i, j}\lambda^-_{k, l}  - \lambda^+_{k, l}\lambda^-_{i, j}
    \end{align}
    Since either $\lambda^+_{i, j} = 0$ or $\lambda^-_{i, j} = 0$, the constraint $\lambda^+_{i, j}\lambda^-_{k, l}  + \lambda^+_{k, l}\lambda^-_{i, j} \geq 0$ can be decomposed into two separate constraints $\lambda^+_{i, j}\lambda^-_{k, l} \geq 0$  and $\lambda^+_{k, l}\lambda^-_{i, j} \geq 0$.
    
    In summary, the equivalent atomic constraints are:
    \begin{align} 
        &\lambda^+_{i, j} \geq 0, \quad &&\lambda^-_{i, j} \geq 0, \quad \lambda^+_{i, j} \lambda^-_{i, j} = 0 \quad  \forall (i, j)\\
        &\lambda^+_{i, j}\lambda^-_{k, l} \geq 0, \quad &&\lambda^-_{i, j} \lambda^+_{k, l} \geq 0\quad \forall (i, j) \neq (k, l).
    \end{align}
    
    \subsubsection{Extensions to \cite{FazlyabMorariEtAl2022}}
    
    \paragraph{Direct Extensions.}
    \cite{NewtonPapachristodoulou2021} showed that the formulation proposed in \cite{FazlyabMorariEtAl2022} exhibits chordal sparsity. They exploit this insight by using an off-the-shelf solver, SparseCoLo, capable of taking advantage of chordal sparsity. This approach offers an advantage for deep networks. Because the formulation is the same as \cite{FazlyabMorariEtAl2022}, it exhibits the same relaxation gap.

\subsection{Rank and Exactness}\label{subsec:rank}
    \cite{Zhang2020} defined an SDP relaxation to be tight if it has a unique rank-one solution. The rationale for such a definition is because if an interior-point method converges to a maximum rank solution, the optimal solution will not be rank-one unless it is unique. 
    
    Based on Theorem \ref{thm:verif_cpp}, however, we propose a less restrictive definition of tightness for SDP-based relaxations. Specifically, we say that an SDP-relaxation is tight if the optimal value is equal to that of \eqref{opt:linear_exact}.
    This definition is motivated by the fact that Theorem \ref{thm:verif_cpp} does not preclude high rank solutions. As a consequence, if the optimal solution for Problem \eqref{opt:0SOS}(0-SOS) is feasible for Problem \eqref{opt:CPP} (CPP), then it is exact, regardless of its rank.

    In this section, we show that the revised definition is not vacuous, and there are in fact verification instances where the optimal solution of the SDP-relaxation is exact despite not being rank-one, and show how the optimal solutions of \eqref{opt:linear_exact} can be recovered from the SDP solution.

    Consider the ReLU network given by the following weights and biases,
    \begin{align}
        W[0] &= 
        \begin{pmatrix} 1 & 1\\ 1 & -1 \\ -1 & 1 \\ -1 & -1 \end{pmatrix}, \quad 
        && b[0]  = \begin{pmatrix} 0 \\ 0 \\ 0 \\ 0 \end{pmatrix} \label{eq:ex_nw_start}\\
        W[1] &= 
        \begin{pmatrix}
        1 & 0 & 0 & 0 \\ -1 & 1 & 0 & 0\\ 0 & 0 & 1 & 0\\ 0 & 0 & -1 & 1 
        \end{pmatrix}, 
        \quad 
        && b[1] =
        \begin{pmatrix} 2.1 \\ 0 \\ 2.1 \\ 0\end{pmatrix}\\
        W[2] &= 
        \begin{pmatrix}
        1 & 1 & 0 & 0\\ -1 & -1 & 1 & 1
        \end{pmatrix}
        , 
        \quad 
        && b[2] = \begin{pmatrix} 0 \\ 0 \end{pmatrix} \\
        W[3] &= 
        \begin{pmatrix}
        -1 \\ -1
        \end{pmatrix}, \quad 
        && b[3] = \begin{pmatrix} 2.1 \end{pmatrix} \label{eq:ex_nw_end}
    \end{align}
    and the safety rule to be verified defined by:
    \begin{align}
        \mathcal{X} &:= \{z_{0, *} \mid -1 \leq z_{0, *} \leq 1\} \subset \mathbb{R}^2\\
        \mathcal{Y} &:= \{z_{3, *} \mid z_{3, *} \geq -2.1\} \subset \mathbb{R}.
    \end{align}
    
    We form the 0-SOS relaxation, \eqref{opt:0SOS}, and solve the SDP using SCS; this had an optimal value of $-2.0002017$.
    
\begin{figure}[ht]
        \centering
        \includegraphics[width=0.5\textwidth]{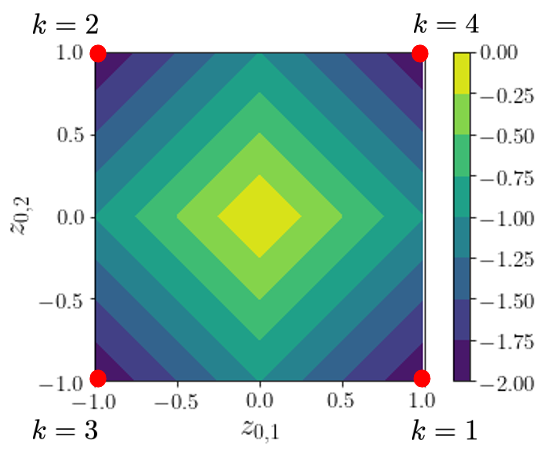}
        \caption{This figure plots the level set of the network defined by the weights and biases in Equations \eqref{eq:ex_nw_start}-\eqref{eq:ex_nw_end}, as well as the inputs, \eqref{eq:sol1}-\eqref{eq:sol4}, corresponding to individual factors in the SDP solution.}
        \label{fig:heatmap}
    \end{figure}

    We then factored the optimal solution using the using the Alternate Least Square Using Projected Gradient Descent algorithm---\cite{Lin2007} (implemented in the NMF.jl package). We found a rank-four factorization with a maximum entry-wise error of $0.0001158$, and objective value of $-2.0002022$. The corresponding values of $\lambda^{\pm, (k)}_{0, *}$ and $\xi^{(k)}$ are as follows:
    \begin{align}
        &k = 1 \quad  
        &&\lambda^{+,(1)}_{0,*} = 
        \begin{pmatrix} 0.491844 \\ 3.32366\times 10^{-5} \end{pmatrix}, \quad 
        &&&\lambda^{-,(1)}_{0,*} = 
        \begin{pmatrix} 0.0  \\ 0.491843 \end{pmatrix}\quad 
        &&&& \xi^{(1)} = 0.491849 \label{eq:sol1}\\
        &k = 2 \quad          
        &&\lambda^{+,(2)}_{0,*} = 
        \begin{pmatrix} 4.76418\times10^{-5} \\ 0.504452 \end{pmatrix}, \quad 
        &&&\lambda^{-,(2)}_{0,*} = 
        \begin{pmatrix} 0.504447  \\ 0.0  \end{pmatrix}\quad 
        &&&& \xi^{(2)} = 0.504478\\
        &k = 3 \quad          
        &&\lambda^{+,(3)}_{0,*} = 
        \begin{pmatrix} 0.0 \\ 0.0 \end{pmatrix}, \quad 
        &&&\lambda^{-,(3)}_{0,*} = 
        \begin{pmatrix} 0.522856 \\ 0.522838  \end{pmatrix}\quad 
        &&&& \xi^{(3)} = 0.522831\\
        &k = 4 \quad          
        &&\lambda^{+,(4)}_{0,*} = 
        \begin{pmatrix} 0.411613  \\ 0.411621  \end{pmatrix}, \quad 
        &&&\lambda^{-,(3)}_{0,*} = 
        \begin{pmatrix} 0.0 \\ 0.0  \end{pmatrix}\quad 
        &&&& \xi^{(4)} = 0.411569 \label{eq:sol4}
    \end{align}
    It can be seen that each solution input corresponds to a corner of $\mathcal{X}$. Each of these solutions are plotted alongside the level sets of the network in Figure \ref{fig:heatmap}. In this example, the optimal SDP-relaxation had a non-negative factorization, but this is not generally guaranteed. This does suggest, however, that one should attempt to find a non-negative factorization of the 0-SOS relaxation as a certificate of optimality.

\end{document}